\documentclass[a4paper,10pt]{amsart}
\usepackage{amsmath,amsthm,amssymb,enumerate}
\usepackage{graphicx}

\numberwithin{equation}{section}
\numberwithin{figure}{section}

\newtheorem{thm}{Theorem}[section]
\newtheorem{lem}[thm]{Lemma}
\newtheorem{prop}[thm]{Proposition}
\newtheorem{cor}[thm]{Corollary}
\newtheorem{definition}[thm]{Definition}
\newtheorem{example}[thm]{Example}

\newenvironment{df}{\begin{definition}\rm}{\end{definition}}

\newcommand{\Proof}{\noindent{\it Proof.} }
\newcommand{\statement}{\noindent{\bfseries Remark.} }
\newcommand{\acknowledgments}{\noindent{\bfseries Acknowledgments.} }

\title{SURFACE LINKS WHICH ARE COVERINGS 
OVER \\
THE STANDARD TORUS}

\author{Inasa Nakamura$^*$}
\thanks{$^*$Supported by GCOE, Kyoto University}

\date{}



\address{Research Institute for Mathematical Sciences, Kyoto University, Kyoto 606-8502, 
JAPAN}
\email{inasa@kurims.kyoto-u.ac.jp}

\begin{document}

\begin{abstract}
We introduce a new construction of a surface link in the 4-space. 
We construct a surface link as a branched covering 
over the standard torus, which we call a torus-covering link. 
We show that 
a certain torus-covering $T^2$-link is equivalent to 
the split union of 
spun $T^2$-links and turned spun $T^2$-links. 
We show that a certain torus-covering $T^2$-link has a non-classical link group. 
We give a certain class of ribbon torus-covering $T^2$-links. 
We present the quandle cocycle invariant of a certain torus-covering $T^2$-link 
obtained from a classical braid, 
by using the quandle cocycle invariants of the closure of the braid. 
\end{abstract}

\maketitle
%

\section{Introduction}
%
%
A {\it surface link} is the image of a smooth embedding of a closed surface into 
the Euclidean 4-space $\mathbb{R}^4$. 
It is known \cite{Kamada1, Kamada3} that any oriented surface link can be presented by the 
closure 
of a surface braid. 
Here, the closure of a surface braid is a surface link of the following form.
Let $S^2$ be a standard 2-sphere in $\mathbb{R}^4$, i.e. the boundary of a standard 3-ball 
in $\mathbb{R}^3 \times \{0\}$. 
The {\it closure of a surface braid} is a surface link embedded in a tubular 
neighborhood $N(S^2)$ of $S^2$ in 
such a way that the projection of it to $S^2$ is a branched covering over $S^2$. 
We identify $N(S^2)$ with $I \times I \times S^2$, where 
$I$ is an interval. 
For a surface link $S$ of such a form,
we consider the singular set of the image of $S$ by the projection to $I \times S^2$, and 
the image of this singular set by the projection to $S^2$ forms a graph on $S^2$. 
An {\it $m$-chart on} $S^2$ is such a graph with certain additional data. 
We can present the original surface link by its $m$-chart on $S^2$ 
(\cite{Kamada2, Kamada3}).

 In this paper we introduce a \lq\lq torus-covering link" as 
a new construction of a surface link, 
by considering a standard torus instead of a standard 2-sphere. 
Let $T$ be a standard torus in $\mathbb{R}^4$, i.e. 
the boundary of a standard solid torus in $\mathbb{R}^3 \times \{0\}$. 
A {\it torus-covering link} is a surface link embedded in a tubular 
neighborhood $N(T)$ of $T$ in 
such a way that the projection of it to $T$ is a branched covering over $T$. 
For a surface link of such a form,
we can define its {\it $m$-chart on} $T$ in the same way as above. 
 A torus-covering link can be presented by an $m$-chart on $T$. 
The aim of this paper is to study various aspects of torus-covering links. 
 
 We introduce an equivalence relation, called the $t$-equivalence, 
among $m$-charts on $T$, 
and show that two torus-covering links are equivalent if their 
$m$-charts on $T$ are $t$-equivalent (Theorem \ref{0706-2}). 
A {\it $T^2$-link} is a surface link whose components are homeomorphic 
to tori. 
We show that a torus-covering $T^2$-link is determined from two commutative classical $m$-braids 
(Lemma \ref{55}), which we call {\it basis $m$-braids}, and 
we denote by $\mathcal{S}_m(a,b)$ the torus-covering $T^2$-link with 
basis $m$-braids $a$ and $b$. 
A vertex of degree one (resp. six) of an $m$-chart is called 
a {\it black vertex} (resp. {\it white vertex}). 
A torus-covering $T^2$-link is presented by an $m$-chart on $T$ without black vertices 
(Lemma \ref{56}). 
We show that 
 an $m$-chart on $T$ with neither black nor white vertices presents the split union 
of spun $T^2$-links and turned spun $T^2$-links (Theorem \ref{0706-1}). 

The {\it link group} of a surface link or a classical link is the fundamental 
group of the link exterior. 
First we calculate the link group of $\mathcal{S}_m(a,b)$ (Proposition \ref{Lem4-1}). 
It is known \cite{Livingston, Boyle} that 
a spun $T^2$-link or a turned spun $T^2$-link has a classical link group; 
thus the split union of spun $T^2$-links and turned spun $T^2$-links
also has a classical link group. 
We will show that a certain 2-component torus-covering $T^2$-link has a 
non-classical link group  
(Theorem \ref{Thm4-6}). 
We show its knot version as well: a certain torus-covering $T^2$-knot 
has a non-classical knot group (Theorem \ref{Thm4-9}). 
As a corollary, we can see that the torus-covering $T^2$-link of 
Theorems \ref{Thm4-6} or \ref{Thm4-9}
   is not equivalent to the split union of 
 spun $T^2$-links and turned spun $T^2$-links (Theorem \ref{0730-1}). 

An oriented surface link is called {\it ribbon} 
if it is the boundary of an immersed 3-manifold with \lq\lq ribbon singularities" (\cite{Yanagawa}). 
We give a certain class of ribbon torus-covering $T^2$-links 
(Theorem \ref{Prop2-10}). 
As a corollary, we can see that 
the torus-covering $T^2$-link of Theorem \ref{0730-1} 
   is ribbon (Corollary \ref{0704-1}).

It is known \cite{Asami-Satoh} that the quandle cocycle invariant
of a twist-spun 2-knot of a classical knot $K$ 
can be presented by using the quandle cocycle invariants of a 1-tangle whose closure is $K$. 
 From a similar viewpoint, we expect that an invariant of $\mathcal{S}_m(b, \Delta^{2n})$ 
can be presented by using invariants of an $m$-braid $b$, 
where $\Delta$ is a half twist of a bundle of $m$ 
parallel strands. 
In Theorem \ref{0602-t} we present the quandle cocycle invariant 
of $\mathcal{S}_m(b, \Delta^{2n})$, 
by using the quandle cocycle invariants of the closure of $b$. 
In Theorem \ref{Thm2-11}, we calculate some concrete examples of Theorem \ref{0602-t}. 
  They give torus-covering $T^2$-knots whose triple point 
 numbers are positive (Corollary \ref{0526}).

The paper is organized as follows. 
In Section \ref{1}, we define a torus-covering link (Definition \ref{Def2-1}) 
and show Theorem \ref{0706-2}. 
Further we 
study torus-covering $T^2$-links and show Theorem \ref{0706-1}. 
In Section \ref{KnotGroup}, 
we study link groups of torus-covering $T^2$-links. 
We show Theorems \ref{Thm4-6} and \ref{Thm4-9}. 
Further, we show 
Theorem \ref{0730-1}. 
In Section \ref{ribbon}, we show Theorem \ref{Prop2-10}. 
In Section \ref{TriplePoint}, we calculate the quandle cocycle invariants and 
show Theorem \ref{0602-t}. 
Further we show 
Theorem \ref{Thm2-11}, by using Mochizuki's 3-cocycle. 
 
%
%
\section{Torus-covering links} \label{1}
%
A braided surface over a 2-disk was defined in \cite{Rudolph, Kamada3}. 
A surface braid is a braided surface with some boundary condition, and 
a notion of an $m$-chart on a 2-disk was introduced \cite{Kamada92, Kamada3} to present 
a simple surface braid. 
Equivalent simple surface braids have distinct chart presentations. 
The notion of C-move equivalence between two $m$-charts on a 2-disk was introduced 
\cite{Kamada92, Kamada2, Kamada3} to 
give the equivalence class of an $m$-chart which represents 
the equivalence class of a simple surface braid. 
In this section, we modify the definitions to define a braided surface $S$ 
over a closed surface 
$\Sigma$, 
an $m$-chart on $\Sigma$ which presents $S$, and 
the notion of C-move equivalence between two $m$-charts on $\Sigma$. 
Using these terms, 
we define a torus-covering link, which is presented by an $m$-chart on the standard torus. 
We define $t$-equivalence between two $m$-charts, and show that 
the torus-covering links are equivalent 
if $m$-charts of them are $t$-equivalent (Theorem \ref{0706-2}). 
Further we 
study torus-covering $T^2$-links. 
A torus-covering $T^2$-link is presented by an $m$-chart on $T$ without black vertices 
(Lemma \ref{56}). 
We show Theorem \ref{0706-1}: 
an $m$-chart on $T$ with neither black nor white vertices presents the split union 
of spun $T^2$-links and turned spun $T^2$-links. 

We work in the smooth category, and we assume that embeddings are locally flat. 
Let $D^2=I \times I$, where $I=[0,1]$. 
A {\it surface link} is the image of a smooth embedding of a closed surface into 
$\mathbb{R}^4$. 
Two surface links are said to be {\it equivalent} if one is taken to the other 
by an orientation-preserving self-diffeomorphism of $\mathbb{R}^4$. 

\begin{df}
A closed surface $S$ embedded in $D^2 \times \Sigma$ 
is called a {\it braided surface over 
$\Sigma$} 
of degree $m$
if 
 $p_{\Sigma} |_{S} \,:\, S \rightarrow \Sigma$ is a branched covering map of 
degree $m$, where 
$p_{\Sigma} \,:\, D^2 \times \Sigma 
\rightarrow \Sigma$ is 
the projection to the second factor. 
A braided surface $S$ is called \textit{simple} 
if $\#(S \cap p_{\Sigma}^{-1}(x))=m-1$ or $m$ for each $x \in \Sigma$. 
Take a base point $x_0$ of $\Sigma$. 
Two braided surfaces over $\Sigma$ of degree $m$ are 
{\it equivalent} if there is a fiber-preserving ambient isotopy 
of 
$D^2 \times \Sigma$ rel $p_{\Sigma}^{-1}(x_0)$ 
which carries one to the other. 
\end{df}

When a simple braided surface $S$ is given, we obtain a graph on $\Sigma$, as follows. 
Consider the singular set $\mathrm{Sing}(p_1(S))$ of the image of $S$ by 
the projection $p_1$ to $I \times \Sigma$. 
Perturbing $S$ if necessary, 
we can assume that 
$\mathrm{Sing}(p_1(S))$ consists of 
double point curves, 
triple points, and branch points. 
Moreover we can assume that the singular set of 
the image of $\mathrm{Sing}(p_1(S))$ 
by the projection to 
$\Sigma$ consists of a finite number of double points such that the preimages 
belong to double point curves of $\mathrm{Sing}(p_1(S))$. 
Thus 
the image of $\mathrm{Sing}(p_1(S))$ by the projection to 
$\Sigma$ forms a finite graph $\Gamma$ on $\Sigma$ such that 
the degree of its vertex is either $1$, $4$ or $6$. 
An edge of $\Gamma$ corresponds 
to a double point curve, and a vertex of degree $1$ (resp. $6$) 
corresponds to a branch point (resp. triple point). 

For such a graph $\Gamma$ obtained from a simple braided surface $S$, 
we give orientations and labels to the edges of $\Gamma$, as follows. 
Let us consider a path $l$ in $\Sigma$ such that 
$l \cap \Gamma$ is a point $P$ of an edge $e$ of $\Gamma$. 
Then $S \cap p_{\Sigma}^{-1} (l)$ is a classical $m$-braid with one crossing 
in $p_{\Sigma}^{-1}(l)$ 
such that $P$ corresponds to the crossing of the $m$-braid. 
Let $\sigma_1, \sigma_2, \ldots, \sigma_{m-1}$ 
be the standard generators of the $m$-braid group $B_m$. 
Let $\sigma_{i}^{\epsilon}$ 
($i \in \{1,2,\ldots, m-1\}$, 
$\epsilon \in \{+1, -1\}$) be 
the presentation of $S \cap p_{\Sigma}^{-1}(l)$. 
Then label the edge $e$ by $i$, and moreover give $e$ an orientation such that 
the normal vector of $l$ corresponds (resp. 
does not correspond) 
to the orientation of $e$ if $\epsilon=+1$ (resp. $-1$). 
We call such an oriented and labeled graph an {\it $m$-chart of $S$}. 
   \\
 
 In general, we define an $m$-chart on $\Sigma$ as follows. 

\begin{figure}
 \includegraphics*{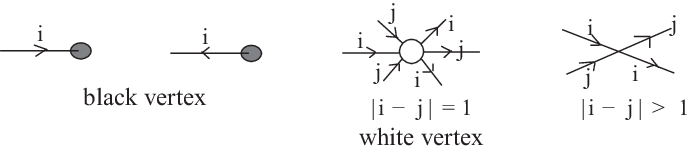}
 \caption{Vertices in an $m$-chart}
 \label{Fig1-1}
 \end{figure}

\begin{df}
Let $m$ be a positive integer, and let 
$\Gamma$ be a finite graph on $\Sigma$.
Then $\Gamma$ is called an {\it m-chart on $\Sigma$} if 
it satisfies the following conditions: 

\begin{enumerate}[(i)]
 \item Every edge is oriented and labeled by an element of 
       $\{1,2, \ldots, m-1\}$. 
 \item Every vertex has degree $1$, $4$, or $6$.
 \item  The adjacent edges around each vertex are oriented and labeled as shown in 
Fig. \ref{Fig1-1}, 
where we depict a vertex of degree 1 by a black vertex, and a vertex of degree 
6 by a white vertex. 
 \end{enumerate}
 \end{df}
 
 When an $m$-chart $\Gamma$ on $\Sigma$ is given, 
 we can reconstruct a simple braided surface $S$ over $\Sigma$ 
as follows. 
Let $N(\Gamma)$ be a neighborhood of $\Gamma$ in $\Sigma$. 
Let us consider a trivial braided surface $S=Q_m \times (\Sigma -N(\Gamma))$ 
over $\Sigma -N(\Gamma)$, where $Q_m$ is a set of $m$ interior points of $D^2$. 
We extend $S$ over a neighborhood of each edge as follows. 
Identify a neighborhood of an edge $e$ with $I \times I$ such that 
 $e$ is identified with $\{1/2\} \times I$. 
Let $i$ be the label attached to $e$, and let $\epsilon=+1$ (resp. $-1$) 
if the orientation of $e$ corresponds (resp. does not correspond) 
to the orientation of $\{0\} \times I$. 
Then let the braided surface $S$ over the neighborhood of $e$ be 
the braided surface which has a presentation 
$\sigma_i^\epsilon \times I$ and 
the image of the double point curve of $p_1(S)$ by the projection to $\Sigma$ 
is $e$. 
Since $\Gamma$ is as in Fig. \ref{Fig1-1} around each vertex, 
$S$ can be extended naturally over a neighborhood 
of each vertex. See \cite{Carter-Saito, Kamada92-2, Kamada3} 
for more details. 
Thus we can construct a simple braided surface $S$ over $\Sigma$ 
such that the original $m$-chart $\Gamma$ is an $m$-chart of $S$. 
\\

\begin{figure}
 \includegraphics*{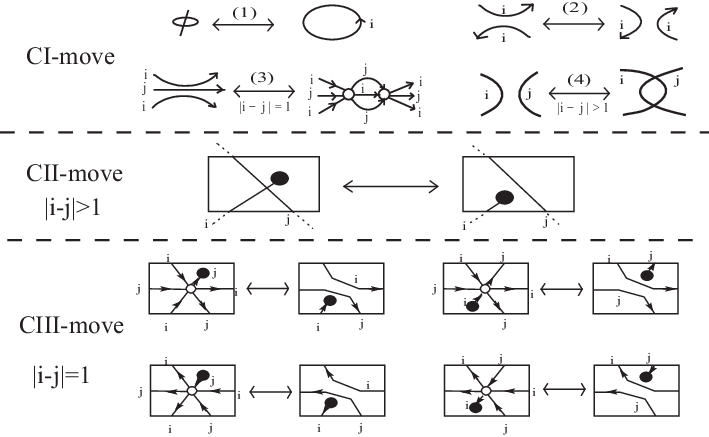}
\caption{CI, CII, CIII-moves. For CI-moves, we give only several examples.} 
\label{cmove}
 \end{figure}

Two $m$-charts on $\Sigma$ are {\it C-move equivalent} if 
they are related by a finite sequence of 
ambient isotopies of $\Sigma$ rel $p_{\Sigma}^{-1}(x_0)$ and CI, CII, CIII-moves shown in Fig. \ref{cmove}; 
see \cite{Kamada3} for the complete set of CI-moves. 
 It is shown as a minor modification of \cite{Kamada92, Kamada2, Kamada3} that 
two simple braided surfaces over $\Sigma$ of degree $m$ are equivalent if and only if 
$m$-charts of them are C-move equivalent. 

Now we define torus-covering links. 
 Let $T$ be 
the standard torus in $\mathbb{R}^4$, i.e. 
the boundary of the standard solid torus in $\mathbb{R}^3 \times \{0\}$. 
Let us fix a point $x_0$ of $T$, and take a meridian 
$\mathbf{m}$ and a longitude $\mathbf{l}$ of $T$ with the base point $x_0$. 
A meridian is an oriented simple closed curve 
on $T$ which bounds the 2-disk of the solid torus whose boundary is $T$. 
A longitude is an oriented simple closed curve 
on $T$ which is null-homologous 
in the complement of the solid torus in the three space $\mathbb{R}^3 \times \{0\}$. 
%
Let 
$N(T)$ be a 
tubular neighborhood of $T$ in $\mathbb{R}^4$. 
Since $T$ is the boundary of the standard solid torus in $\mathbb{R}^3 \times \{0\}$, 
the normal bundle of $T$ in $\mathbb{R}^3 \times \{0\}$ is a trivial bundle. 
Let us 
identify it with $I \times T$. 
Then let us identify 
$N(T)$ with $I \times I \times T$, where 
the second $I$ is an interval in the fourth axis of $\mathbb{R}^4$. 
 From now on, we identify $N(T)$ with $D^2 \times T$. 

\begin{df} \label{Def2-1} 
A {\it torus-covering link} is a surface link in $\mathbb{R}^4$ 
presented by a simple braided surface over $T$, where we regard the braided surface as in 
$N(T) \subset \mathbb{R}^4$. 
 \end{df}
  As we mentioned, for two $m$-charts on $T$, 
their presenting braided surfaces over $T$ are equivalent if 
the $m$-charts are C-move equivalent.  
Hence it follows that 
for two $m$-charts on $T$, 
their presenting torus-covering links are equivalent if 
the $m$-charts are C-move equivalent. 
 Since each component of a torus-covering link is a branched cover over a torus $T$, 
 each component of a torus-covering link 
 is of genus at least one. 
See Propositions \ref{eg}, \ref{2eg} and \ref{3eg} for some examples of torus-covering links. 
Note that it is known \cite{Berstein-Edmonds, B-E84} that any 
braided surface over $T$ is approximated by 
a simple braided surface over $T$. 

Regarding $S^4$ as the one-point compactification of $\mathbb{R}^4$, 
we regard a surface link as in $S^4$. 
Then $N(T)=D^2 \times T$ is embedded in 
 $S^4$. 
Let $\mathbf{r}=\partial D^2 \times \{0\} \times \{0\}$ 
 be a 
curve on $\partial N(T)=\partial D^2 \times T$. 
Put $E^4={\rm cl} (S^4-N(T))$. 
Let $r$, $m$ and $l$ be the curves on $\partial E^4$, 
which are identified with $\mathbf{r}$, $\mathbf{m}$ and $\mathbf{l}$ under 
the natural identification map 
$i\,:\,\partial N(T) \rightarrow \partial E^4$. 
The curves $\mathbf{r}$, $\mathbf{m}$ and $\mathbf{l}$ represent a basis of $H_1(\partial E^4 \,;\, \mathbb{Z})$. 
Let $f \,:\, \partial E^4 \rightarrow E^4$ be a diffeomorphism such that 
$f_* 
(
\begin{array}{ccc}
 \mathbf{r} & \mathbf{m} & \mathbf{l}
\end{array}
)
=
(
\begin{array}{ccc}
 \mathbf{r} & \mathbf{m} & \mathbf{l}
\end{array}
)
A^f$, 
where $A^f \in GL(3,\mathbb{Z}) \cong \pi_0 \mathrm{Diffeo}(\partial E^4)$. 
It is known \cite{Montesinos} that 
the map $f$ can be extended to a self-diffeomorphism of $E^4$ 
if and only if $A^f \in H$, where
\[
H=\left\{
\left(
\begin{array}{ccc}
\pm 1 & 0 & 0 \\
* & \alpha & \gamma \\
* & \beta & \delta
\end{array}
\right)
\in GL(3, \mathbb{Z}) \, ; \, \alpha+\beta+\gamma+\delta \equiv 0
\ \ (\bmod \ 2)
 \right\}. 
 \] 
Using this fact, we introduce an equivalence relation 
between two $m$-charts on $T$.

\begin{df} \label{te}
We say that two $m$-charts on $T$ are {\it t-equivalent} if they are related by 
a finite sequence of ambient isotopies of $T$, C-moves and 
a self-diffeomorphism of $T$ given by an element of 
\[\left\{
\left(
\begin{array}{cc}
\alpha & \gamma \\
\beta & \delta
\end{array}
\right)
\in GL_+(2, \mathbb{Z}) \, ; \, \alpha+\beta+\gamma+\delta \equiv 0
\pmod{2}
 \right\}. 
\]
\end{df}

\begin{thm} \label{0706-2}
Two torus-covering links in $S^4$ 
are equivalent if $m$-charts of them are t-equivalent. 
\end{thm}

\Proof
Since C-move equivalent $m$-charts present equivalent torus-covering links, 
it suffies to show in the cases for an ambient isotopy of $T$ and 
a self-diffeomorphism $g$ of $T$ of Definition \ref{te}. 

An ambient isotopy of $T$ induces a fiber-preserving ambient isotopy of $N(T)$ 
which relates the torus-covering links. 
This can be extended to the whole space by the Isotopy Extension Theorem (see \cite{Hirsch}). 
The terminal map of the resulting isotopy is an orientation-preserving diffeomorphism; thus 
the torus-covering links are equivalent.  

Let $\Gamma$ be an $m$-chart on $T$. 
Let $S$ and $S^\prime$ be the torus-covering links presented by 
$\Gamma$ and $g(\Gamma)$ respectively. 
Let $g^\prime$ be a self-diffeomorphism of $N(T)$ induced by $g$, i.e. 
$g^\prime=\mathrm{id}_{D^2} \times g \,:\, N(T) \rightarrow 
 N(T)$, where we regard $N(T)$ as $D^2 \times T$. 
Since 
$A^{g^\prime}$ 
is in $GL_+(3,\mathbb{Z})$, 
the map $g^\prime|_{\partial N(T)}$ 
 can be considered as an orientation-preserving self-diffeomorphism of $\partial E^4$. 
Since 
$A^{g^\prime}$ is an element of $H$, 
 $g^\prime|_{\partial N(T)}$ can be extended to $E^4$, and hence to 
$S^4=N(T)\cup _{\partial N(T)} E^4$. 
 This is an orientation-preserving self-diffeomorphism of $S^4$ which maps 
 $S$ to $S^\prime$, and hence $S$ and $S^\prime$ are equivalent in $S^4$. 
\qed
\\

 In particular, we have the following corollary. 
Let $\rho$ (resp. $\tau$) 
be a self-diffeomorphism of $T$ given by 
$\left(
\begin{array}{cc}
 0 & -1 \\
 1 & 0
\end{array}
\right)$
(resp.  $\left(
\begin{array}{cc}
 1 & -1 \\
 0 & 1
\end{array}
\right)$ ). 

\begin{figure}
 \includegraphics*{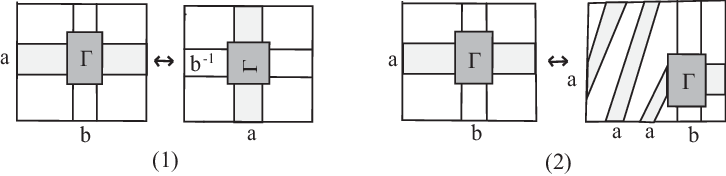}
\caption{Equivalent moves of $m$-charts} 
\label{0629-1}
 \end{figure}

\begin{cor} \label{Prop2-3}
 Two $m$-charts on $T$ present equivalent torus-covering links if they are related by 
a finite sequence of ambient isotopies of $T$, C-moves, and moves 
as in Fig. \ref{0629-1}. 
\end{cor}

\Proof
The moves as in Fig. \ref{0629-1} (1) and (2) are 
related by 
$\rho$ and 
$\tau^2$ respectively. 
They 
give $t$-equivalence between two $m$-charts on $T$. 
Thus 
the conclusion follows from 
Theorem \ref{0706-2}. 
\qed
\\

\statement
 Teragaito \cite{Teragaito} proved the same fact of the above corollary 
for the symmetry-spun 
version. 
The case of $\tau^2$ for turned spun $T^2$-links 
was shown by Boyle \cite{Boyle}, using the result of 
Gluck \cite{Gluck}. 
\\

A {\it torus-covering $T^2$-link} is a torus-covering link whose 
components are homeomorphic to tori. 
\begin{lem} \label{56}
A torus-covering $T^2$-link $S$ is presented by 
 an $m$-chart without black vertices. 
Thus $S$ is an unbranched cover over $T$. 
\end{lem}

\Proof
Let $\Gamma$ be an $m$-chart on $T$ which presents $S$, and 
let $b(\Gamma)$ be the number of black vertices in $\Gamma$. 
Then $S$ is a branched cover over $T$ with $b(\Gamma)$ branch points, and 
 the Euler characteristic of $S$ is equal to $-b(\Gamma)$. 
Further, since the Euler characteristic of tori is equal to zero, 
 we have $b(\Gamma)=0$. 
\qed
 \\

Let us consider a torus-covering $T^2$-link $S$. 
 The intersections 
$S \cap p_T^{-1}(\mathbf{m})$ and 
$S \cap p_T^{-1}(\mathbf{l})$ are 
closures of classical braids. 
Cutting open the solid tori at the 2-disk $p_T^{-1}(x_0)$, 
we obtain a pair of classical braids. 
We call them 
{\it basis braids}. 

\begin{lem} \label{55}
\noindent
\begin{enumerate}
\item
The basis braids of a torus-covering $T^2$-link are commutative. 
\item 
For any commutative $m$-braids $a$ and $b$, 
there exists a unique torus-covering $T^2$-link with basis braids $a$ and $b$. 
\end{enumerate}
\end{lem}

For commutative $m$-braids $a$ and $b$, 
we denote by $\mathcal{S}_m(a,b)$ the torus-covering $T^2$-link with 
basis $m$-braids $a$ and $b$. 
\\

\Proof 
 (1) Let $X_m$ be the configuration space of unordered $m$ distinct points of 
a 2-disk $D^2$, i.e. the set of $m$-element subsets of $D^2$ such that each $m$-element subset 
consists of $m$ distinct points. 
It is known \cite{Artin} (see also \cite{Birman}) %
that $\pi_1(X_m)=B_m$. 
Since a torus-covering $T^2$-link is presented by 
an unbranched covering 
of degree $m$ over $T$ by Lemma \ref{56}, it is presented by 
a map $f\,:\,T \rightarrow X_m$. The induced map $f_*$ gives a homomorphism 
$\mathbb{Z} \oplus \mathbb{Z} \cong \pi_1(T) \rightarrow \pi_1(X_m)=B_m$. 
 Since the basis braids are the images of generators of $\mathbb{Z} \oplus \mathbb{Z}$ 
by this homomorphism, they are commutative.

(2) 
For any commutative $m$-braids $a$ and $b$, let us consider 
a map $f\,:\, \mathbf{m} \cup \mathbf{l} \rightarrow X_m$ 
such that the closed paths $f|_{\mathbf{m}}$ and $f|_{\mathbf{l}}$ in $X_m$ 
represent $a$ and $b$ respectively. 
Since $a$ and $b$ are commutative, $a b a^{-1} b^{-1}$ is isotopic to 
the trivial braid, and it follows that the closed path $l$ in $X_m$ representing 
$a b a^{-1} b^{-1}$ is null-homotopic. 
Hence we can take a 2-disk in $X_m$ such that the boundary is $l$:  
thus $f$ can be extended to a map from $T$ to $X_m$. 
Since $\pi_2(X_m)=0$ (\cite{Fadell-Neuwirth}), such an extension is 
unique (up to equivalence). This means that 
there exists a unique torus-covering $T^2$-link with basis braids $a$ and $b$. 
\qed
\\

By Corollary \ref{Prop2-3}, we have the following corollary. 
\begin{cor} \label{0629-2}
For commutative $m$-braids $a$ and $b$, the following equivalent relations hold:
\begin{equation*} 
\mathcal{S}_m(a,b) \sim \mathcal{S}_m(b^{-1},a), \  \mathcal{S}_m(a,b) \sim \mathcal{S}_m(a, a^2 b). 
\end{equation*}
\end{cor}
\Proof
Let $\Gamma$ be an $m$-chart on $T$ 
without black vertices and with  
basis braids $a$ and $b$. By Lemmas \ref{56} and \ref{55} (2), $\Gamma$ presents $\mathcal{S}_m(a,b)$. 
By Corollary \ref{Prop2-3}, the torus-covering $T^2$-links presented by 
$\Gamma$, $\rho(\Gamma)$ and $\tau^2(\Gamma)$ are 
equivalent. 
The basis braids of $\rho(\Gamma)$ are $b^{-1}$ and $a$. 
Further, 
 the basis braids of $\tau^2(\Gamma)$ are $a$ and $a^2 b$; see Fig. \ref{0629-1}. 
Thus we have the required equivalent relations. 
\qed
\\

Using this corollary, we can for example show the following proposition. 
We remark the result before the statement of Corollary \ref{Cor2-12}. 
An oriented surface link $S$ is {\it invertible} 
if $S$ is equivalent to its orientation-reversed image $-S$. 
\begin{prop} \label{0803-1}
 For any integers $p$ and $n$,  
$\mathcal{S}_4(\sigma_1 \sigma_2^p \sigma_3, \Delta^{2n})$ is invertible.
\end{prop}
\Proof
Let us determine the basis braids of $-\mathcal{S}_m(a,b)$, 
as follows. Put $S=\mathcal{S}_m(a,b)$. 
Then $(-S) \cap p_T^{-1}(-\mathbf{m})$ is the closure of $-a=\bar{a}^{-1}$, and 
$(-S) \cap p_T^{-1}(-\mathbf{l})$ is the closure of $-b=\bar{b}^{-1}$. 
Here 
$-x$ ($x=\mathbf{m}, \mathbf{l}, a, b$) is the orientation-reversed image 
of $x$, and 
$\bar{a}$ (resp. $\bar{b}$) is the $m$-braid obtained from 
$a$ (resp. $b$) by 
replacing $\sigma_i$ with $\sigma_i^{-1}$ for each standard generator $\sigma_i$ 
of $B_m$. 
Hence $(-S) \cap p_T^{-1}(\mathbf{m})$ and 
$(-S) \cap p_T^{-1}(\mathbf{l})$ are the closures of $\bar{a}$ and $\bar{b}$ respectively; 
thus the basis braids of $-S$ are $\bar{a}$ and $\bar{b}$. 
From now on, put $S=\mathcal{S}_4(\sigma_1 \sigma_2^p \sigma_3, \Delta^{2n})$. 
By the above argument and Lemma \ref{55} (2), $-S$ is equivalent to 
$\mathcal{S}_4(\sigma_1^{-1} \sigma_2^{-p} \sigma_3^{-1}, \Delta^{-2n})$. 
Applying the first equivalent relation of Corollary \ref{0629-2} twice, we have 
$\mathcal{S}_m(a,b) \sim \mathcal{S}_m(b^{-1}, a) \sim \mathcal{S}_m(a^{-1}, b^{-1})$. 
 Thus 
$-S \sim \mathcal{S}_4(\sigma_3 \sigma_2^p \sigma_1, \Delta^{2n})$.
Regarding the $i$th string of the basis braids as the $(4-i)$th string ($i=1,\ldots, 4$), 
we can regard the basis braids as $\sigma_1 \sigma_2^p \sigma_3$ and $\Delta^{2n}$ 
respectively; thus $-S \sim \mathcal{S}_4(\sigma_1 \sigma_2^p \sigma_3, \Delta^{2n})=S$. 
  \qed
 \\

We show that some torus-covering $T^2$-links are equivalent to known $T^2$-links. 
Let $b$ an $m$-braid, and let $Q_m$ be the starting point set of $b$. 
Let us denote by $\hat{b}$ the closure of $b$. 
\begin{prop} \label{eg}
The torus-covering $T^2$-link $\mathcal{S}_m(b, e)$ 
 is equivalent to the spun $T^2$-link of 
 $\hat{b}$. 
\end{prop}

The 4-space $\mathbb{R}^4$ is constructed 
by rotating the upper half plane $\mathbb{R}^3_+=\mathbb{R}^2 \times [0, \infty)$ in 
$\mathbb{R}^2 \times \mathbb{R}^2$ around the axis $\mathbb{R}^2 \times \{0\}$. 
This structure is called an {\it open book structure}. 
Let $B^3$ be a 3-ball in $\mathbb{R}^3_+$, and let us naturally identify the orbit of 
$B^3$ with $B^3 \times S^1$ in $\mathbb{R}^4$, where $S^1$ is a circle. 
Let $\pi \,:\, B^3 \times S^1 \rightarrow B^3$ be the projection. 
A surface link $S$ in $B^3 \times S^1$ can be considered as 
$\cup_{t \in S^1} S_t \times \{t\}$, 
where $S_t= \pi(S \cap (B^3 \times \{t\})) \subset B^3$. 
We call the collection $\{S_t\}_{t \in S^1}$ the {\it motion picture} 
of $S$ along $S^1$. 

Let $L$ be a classical link in $B^3$. 
The {\it spun $T^2$-link} of $L$ is the surface link defined by the motion picture 
$S_t=L$ for $t \in S^1$ (\cite{Livingston, Boyle88, Boyle}). 

In an open book structure of $\mathbb{R}^4$, 
 we naturally identify 
the orbit of $B^3$ with $B^3 \times S^1$. 
Let us consider un unknotted circle $S^1 \subset B^3$ with the base point $x_0$. 
The orbit of this $S^1$ is a standard torus. We identify it with $T$, 
by identifying $S^1 \times \{0\} \subset B^3 \times \{0\}$ 
with the meridian $\mathbf{m}$, and $\{x_0\} \times S^1 \subset B^3 \times S^1$ 
with the longitude 
$\mathbf{l}$. 
Further, we identify $N(T)$ with the orbit of the unknotted solid torus 
$N(\mathbf{m})=\pi(p_T^{-1}(\mathbf{m}))$ 
in $B^3$. 
Put $D(x_0)=\pi(p_T^{-1}(x_0))$. 
Let us identify $S^1$ with $[0,1]/\sim$, where $0 \sim 1$. 
Let us assume that the closure $\hat{b}$ is in the solid torus $N(\mathbf{m})$ 
such that the identified corresponding ends are in $D(x_0)$. 
Further we regard a braided surface over $T$ as in $N(T) \subset \mathbb{R}^4$. 
\\

\noindent
{\it Proof of Proposition \ref{eg}}. 
Let $S$ be the surface defined by 
the motion picture 
$S_t =\hat{b}$ for $t \in S^1$. 
By definition, $S$ is the spun $T^2$-link of $\hat{b}$. 
 Since $S$ is a braided surface over $T$ with no branch points, 
$S$ is a torus-covering $T^2$-link. 
Let us determine the basis braids. 
Since $\pi(S \cap p_T^{-1}(\mathbf{m}))=S_0 \cap N(\mathbf{m})=\hat{b}$, 
one basis braid is $b$. 
By definition, we have 
$S \cap p_T^{-1}(\mathbf{l})=\cup_{t \in [0,1]} (S_t \cap D(x_0)) \times \{t\}$. 
Since $S_t \cap D(x_0)=Q_m$ for any $t$, 
$S \cap p_T^{-1}(\mathbf{l})$ is the closure of the trivial $m$-braid $e=Q_m \times [0,1]$. 
Thus $S$ is a torus-covering $T^2$-link with basis $m$-braids $b$ and $e$, 
and it follows from Lemma \ref{55} (2) that $S$ is equivalent to $\mathcal{S}_m(b,e)$. 
\qed 
\\

Let us identify the 3-ball $B^3$ with the unit ball in the xyz-space. 
Let us rotate a classical link $L$ in $B^3$ around the 
$z$-axis once, and identify the resulting link with the original link. 
The orbit of $L$ forms a surface link, called 
the {\it turned spun $T^2$-link} of $L$ (\cite{Livingston, Boyle}).

 \begin{prop} \label{2eg}
 The torus-covering $T^2$-link $\mathcal{S}_m(b, b)$ 
   is equivalent to the turned spun $T^2$-link of $\hat{b}$. 
\end{prop}
 
\Proof
We can assume that the solid torus $N(\mathbf{m}) \subset B^3$ 
is fixed as a set when we rotate it 
around the $z$-axis. 
Let $\{h_u\}_{u \in [0,1]}$ be an isotopy of $B^3$ which describes the 
rotation of the solid torus $N(\mathbf{m})$ 
around the $z$-axis once. 
Let $S$ be the surface defined by 
the motion picture 
$S_t =\cup_{u \in [0,1]} h_{t}(\hat{b})$ for $t \in S^1$. 
By definition, $S$ 
is the 
turned spun $T^2$-link of $\hat{b}$. 
By the same argument with the proof of Proposition \ref{eg}, 
 $S$ is a torus-covering $T^2$-link with $S \cap p_T^{-1}(\mathbf{m})=
\hat{b}$. 
Regarding $N(\mathbf{m})$ as $D^2 \times S^1$ such that $D(x_0)=D^2 \times \{0\}$, 
we have $S_t \cap D(x_0)=h_t(\hat{b}) \cap (D^2 \times \{0\})=
p(\hat{b} \cap (D^2 \times \{t\})) \times \{0\}$ for $t \in [0,1]$, 
where $p \,:\, D^2 \times S^1 \rightarrow D^2$ is the projection; thus 
$S \cap p_T^{-1}(\mathbf{l})=\cup_{t \in [0,1]} (S_t \cap D(x_0)) \times \{t\}=\hat{b}$. 
Thus the basis braids of $S$ are $b$ and $b$, and 
$S \sim \mathcal{S}_m(b,b)$ by Lemma \ref{55} (2). 
\qed
\\

Let $L$ be a classical link in $B^3$ such that 
rotating $L$ around the 
$z$-axis by $2 k \pi /n$ results in the original $L$ as a set. 
Then let us rotate $L$ around the 
$z$-axis by $2k \pi /n$, and identify the resulting link with the original link. 
The orbit of $L$ forms a surface link, 
 called a {\it symmetry-spun $T^2$-link} (\cite{Teragaito}).

Let $b^n$ be an $m$-braid in $D^2 \times [0,1]$ such that 
$b^n \cap (D^2 \times I_j)=b$, where $I_j=[(j-1)/n,\,  j/n]$ 
($j=1,2,\ldots,n$).

\begin{prop} \label{3eg}
  The torus-covering $T^2$-link $\mathcal{S}_m(b^{n}, b^{k})$ 
 is equivalent to a symmetry-spun $T^2$-link, which is constructed by turning $\widehat{b^{n}}$ 
by $2 k \pi /n$ around the axis while spinning. 
 \end{prop}
It is known \cite[Theorem 8]{Teragaito} that the symmetry-spun $T^2$-link 
$\mathcal{S}_m(b^{n}, b^{k})$ is equivalent to either $\mathcal{S}_m(b^{r}, e)$ or 
$\mathcal{S}_m(b^{r}, b^r)$, where $r=\mathrm{gcd}(n,k)$. 
This can be shown by Corollary \ref{0629-2}, too. 
\\

\Proof
Let $\{h_u\}_{u \in [0,1]}$ be an isotopy of $B^3$ which 
describes the rotation of the solid torus $N(\mathbf{m})$ 
around the axis by $2k \pi /n$. 
Let $S$ be the surface defined by 
the motion picture 
$S_t =\cup_{u \in [0,1]} h_{t}(\hat{b})$ for $t \in S^1$. 
By definition, $S$ 
is the symmetry-spun 
 $T^2$-link in question. 
By the same argument with the proof of Proposition \ref{eg}, we can see that 
 $S$ is a torus-covering $T^2$-link with $S \cap p_T^{-1}(\mathbf{m})=
\widehat{b^n}$. 
Using the same notations with the proof of Proposition \ref{2eg}, 
we have $S_t \cap D(x_0)=h_t(\widehat{b^n}) \cap (D^2 \times \{0\})=
p(\widehat{b^n} \cap (D^2 \times \{kt/n\})) \times \{0\}$ for $t \in S^1$. 
Thus 
$S \cap p_T^{-1}(\mathbf{l})=\cup_{t \in S^1} (S_t \cap D(x_0)) \times \{t\}=\widehat{b^k}$. 
Thus the basis braids of $S$ are $b^n$ and $b^k$, and 
$S \sim \mathcal{S}_m(b^n,b^k)$ by Lemma \ref{55} (2). 
\qed
\\

Let us call an edge of an $m$-chart a {\it loop} if it is connected with no vertices. 
Let us consider an $m$-chart on $T$ with no vertices. 
\begin{prop} \label{0628-3}
An $m$-chart on $T$ with no
vertices presents a spun $T^2$-link or a turned spun $T^2$-link. 
\end{prop}
\begin{figure}
 \includegraphics*{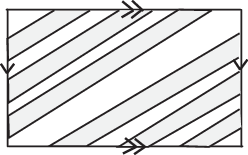}
\caption{} 
\label{0716-1}
 \end{figure} 
\Proof
By Lemma \ref{56}, an $m$-chart with no vertices presents a torus-covering $T^2$-link. 
 Let us determine the basis braids, as follows.  
An $m$-chart with no vertices consists of a finite number of loops. 
By Lemma \ref{0729-1}, we can assume that 
 any loop does not bound a 2-disk in $T$. 
Then, by an ambient isotopy of $T$, we can make all the loops parallel and moreover 
geodesic. 
The $m$-chart is as in Fig. \ref{0716-1}; thus 
the basis braids presented by $\Gamma$ 
are $b^{n}$ and $b^{k}$, for an $m$-braid $b$ and integers $n$ and 
$k$. 
 Its presenting torus-covering $T^2$-link 
is equivalent to a symmetry-spun $T^2$-link by Proposition \ref{3eg}; 
thus it is equivalent to either a spun $T^2$-link or a turned spun $T^2$-link. 
\qed

\begin{lem} \label{0729-1}
An $m$-chart on $T$ with no vertices is C-move equivalent to an $m$-chart 
such that each loop 
does not bound a 2-disk in $T$. 
\end{lem}
\Proof
For a 2-disk $D$ in $T$ such that $\partial D$ intersects an $m$-chart $\Gamma$ 
transversely, if there are no black vertices in $\Gamma \cap D$, then 
by a CI-move 
we can redraw the $m$-chart within $D$ as we like as long as it has no 
black vertices (see \cite{Kamada3}). 
 Hence, if there is a loop which bounds a 2-disk in $T$, then 
we can remove it by applying a CI-move 
around the neighborhood of the loop. 
\qed
\\

Proposition \ref{0628-3} can be extended to 
an $m$-chart on $T$ with neither black nor white vertices, as follows. 
The {\it split union} of two surface links $S_1$ and $S_2$ is a 
surface link presented by the union of the copies of $S_1$ and $S_2$ 
such that for a 3-sphere $S^3$ embedded in $\mathbb{R}^4$, 
$S_1$ is inside of $S^3$ and $S_2$ is outside. 
The 3-sphere $S^3$ is called a {\it separating 3-sphere}. 
\begin{thm} \label{0706-1}
An $m$-chart on $T$ with neither black nor white vertices 
presents either a spun $T^2$-link, a turned spun $T^2$-link, or the split union 
of spun $T^2$-links and turned spun $T^2$-links. 
\end{thm}
Let $\Gamma$ be an $m$-chart with neither black nor white vertices. 
Then every vertex of $\Gamma$ is of degree $4$. 
Since the diagonal edges around a vertex of degree $4$ have the same label and coherent 
orientation (see Fig. \ref{Fig1-1}), 
we can regard 
the union of connected edges of $\Gamma$ with the label $i$ as an oriented immersed circle 
with the label $i$. 
Let us call it just 
an {\it immersed circle}. 
Since the edges around a vertex of degree $4$ have the labels $i$ and $j$ with $|i-j|>1$,   
we can regard $\Gamma$ as 
consisting of immersed circles with transverse intersections such that 
each intersection is formed 
by two immersed circles with the labels $i$ and $j$ with $|i-j|>1$. 
\\

\Proof
Let $\Gamma$ be an $m$-chart with neither black nor white vertices. 
Let $i>0$ be the minimum integer which does not appear as a label of $\Gamma$. 
Let $\Gamma_{<i}$ be a subgraph in $T$ consisting of the edges of $\Gamma$ 
with the labels smaller than $i$. Further, attach to each edge of this $\Gamma_{<i}$ 
the orientation and label induced from 
$\Gamma$. 
Since $\Gamma$ can be regarded as consisting of immersed circles, 
so can 
$\Gamma_{<i}$; 
thus $\Gamma_{<i}$ presents a new $m$-chart on $T$ with neither black nor 
white vertices. 
Similarly, 
let $\Gamma_{>i}$ be a subgraph in $T$ consisting of the edges of $\Gamma$ 
with the labels larger than $i$, with induced labels and orientations. 
 Then $\Gamma_{>i}$ also presents a new $m$-chart on $T$ with neither black nor 
white vertices. 
Since $\Gamma$ has no edge with the label $i$, $\Gamma=\Gamma_{<i} \cup \Gamma_{<i}$. 
By Lemma \ref{0628-1}, 
the torus-covering link presented by $\Gamma$ is equivalent to 
the split union of two torus-covering links presented by $\Gamma_{<i}$ and $\Gamma_{>i}$. 

If there is an immersed circle which bounds a 2-disk in $T$, then we can remove it 
by applying a CI-move 
by the same argument of Lemma \ref{0729-1}. 
 Thus, taking new $i$ if necessary, 
we can assume that $\Gamma_{<i}$ satisfies the conditions of Lemma \ref{0628-2}, i.e. 
(1) any immersed circle does not bound a 2-disk in $T$, and 
(2) there is at least one immersed circle with the label $j$ for every label $j<i$. 
Thus $\Gamma_{<i}$ presents a spun $T^2$-link or a turned spun $T^2$-link 
by Lemma \ref{0628-2}. 
Using induction for $i$, we can see that 
$\Gamma$ presents a spun $T^2$-link, a turned spun $T^2$-link, 
or the split union of spun $T^2$-links and turned spun $T^2$-links. 
\qed

\begin{lem} \label{0628-1}
The torus-covering link presented by $\Gamma_{<i} \cup \Gamma_{>i}$ is equivalent to 
the split union of two torus-covering links presented by $\Gamma_{<i}$ and 
$\Gamma_{>i}$. 
\end{lem}

\begin{figure}
 \includegraphics*{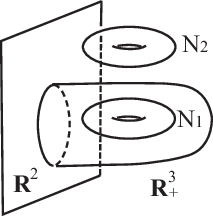}
\caption{} 
\label{0716-2}
 \end{figure} 
\Proof 
Let us denote by $S$, $S_1$, and $S_2$ the 
torus-covering links presented by $\Gamma_{<i} \cup \Gamma_{>i}$, 
$\Gamma_{<i}$, and $\Gamma_{>i}$ respectively. 
Consider the open book decomposition of $\mathbb{R}^4$ as in Proposition \ref{eg}. 
Let $N_1$ (resp. $N_2$) be a solid torus in $N(\mathbf{m})$ which contains 
the $j$th starting point of the basis braids of $S$ for every $j<i$ (resp. $j>i$). 
Since there are no edges of $\Gamma_{<i} \cup \Gamma_{>i}$  with the label $i$, 
we can assume that 
 $S$ is in the orbit 
$(N_1 \cup N_2) \times S^1$ such that $S_1$ and $S_2$ 
are contained in $N_1 \times S^1$ and $N_2 \times S^1$ respectively. 
Let us take a 2-disk in $\mathbb{R}^3_+$ as in Fig. \ref{0716-2}. 
The orbit of this 2-disk forms 
a separating 3-sphere. Thus $S$ is the split union of $S_1$ and $S_2$. 
\qed

\begin{lem} \label{0628-2}
Let $\Gamma$ be an $m$-chart on $T$ with neither black nor white vertices 
satisfying 
(1) any immersed circle does not bound a 2-disk in $T$, and 
(2) there is at least one immersed circle with the label $j$ for every label $j$. 
Then $\Gamma$ presents 
either a spun $T^2$-link or a turned spun $T^2$-link. 
\end{lem}
\Proof
By the definition of an $m$-chart, 
there are no intersections between the immersed circles with the labels $j-1$ and $j$. 
Thus, together with (1), the immersed circles with the labels $j-1$ and $j$ can be made parallel by 
an isotopy of $T$. Hence it follows from (2) that 
each oriented immersed circle goes $n$ times around the meridian $\mathbf{m}$ and 
$k$ times around 
the longitude $\mathbf{l}$, or $-n$ times around $\mathbf{m}$ and 
$-k$ times around $\mathbf{l}$, for some fixed integers $n$ and $k$. 
Since there are only intersections of immersed circles with the labels $i$ and $j$ 
with $|i-j|>1$ as vertices of $\Gamma$, 
we can remove all the vertices by CI-moves of type (4) (see Fig. \ref{cmove}). 
Then the $m$-chart $\Gamma$ presents either a spun $T^2$-link or a turned spun $T^2$-link 
by Proposition \ref{0628-3}. 
\qed
\\

We showed in Theorem \ref{0706-1} that 
an $m$-chart on $T$ with neither black nor white vertices 
presents either a spun $T^2$-link, a turned spun $T^2$-link, or the split union of spun $T^2$-links and turned spun $T^2$-links. 
We will show in Theorem \ref{0730-1} that 
 there is  
a torus-covering $T^2$-link 
which is not 
equivalent to either a spun $T^2$-link, a turned spun $T^2$-link, or 
the split union of spun $T^2$-links and turned spun $T^2$-links. 
Its presenting $m$-chart on $T$ 
does not have black vertices (Lemma \ref{56}) but do have white vertices 
(Corollary \ref{0805-1} (2)).  

\section{Knot Groups and Link Groups} \label{KnotGroup}
 %
From now on throughout this paper, we consider torus-covering $T^2$-links. 
By Lemma \ref{55}, the basis braids of a torus-covering $T^2$-link 
are commutative, and for any commutative $m$-braids $a$ and $b$, 
there exists a unique torus-covering $T^2$-link with basis braids $a$ and $b$. 
In this section, first we compute the link group of $\mathcal{S}_m(a,b)$ (Proposition \ref{Lem4-1}). 
Using this proposition, 
we will show that a certain 2-component torus-covering $T^2$-link  
 has a 
non-classical link group 
(Theorem \ref{Thm4-6}). 
We show its knot version as well: a certain torus-covering $T^2$-knot 
has a non-classical knot group (Theorem \ref{Thm4-9}). 
Further we show that 
the torus-covering $T^2$-link of Theorems \ref{Thm4-6} or \ref{Thm4-9} 
is not equivalent to either a spun $T^2$-link, a turned spun $T^2$-link, or 
the split union of spun $T^2$-links and turned spun $T^2$-links (Theorem \ref{0730-1}). 
 
Before the statement of Proposition \ref{Lem4-1}, we will give the definition of 
Artin's automorphism 
(see \cite{Kamada3}). 
Let $b$ be an 
$m$-braid in a cylinder $D^2 \times [0,1]$, and let $Q_m$ be the starting point set of $b$. 
Let $\{h_u\}_{u \in [0,1]}$ be 
an isotopy of $D^2$ rel $\partial D^2$ 
such that 
$\cup_{u \in [0,1]} h_u (Q_m) \times \{u\}=b$. 
 %
 Let 
$\mathcal{A}^b \,:\, (D^2, Q_m) \rightarrow (D^2, Q_m)$ be the terminal map 
$h_1$, and consider the induced map 
$\mathcal{A}^b_* \,:\, \pi_1(D^2-Q_m) \rightarrow \pi_1(D^2-Q_m)$. 
It is known \cite{Artin} that 
$\mathcal{A}^b$ is uniquely determined from $b$. 
We call $\mathcal{A}^b_*$ {\it Artin's automorphism} associated with $b$. 
Note that $\pi_1(D^2-Q_m)$ is naturally isomorphic to the free group $F_m$ 
generated by the standard generators $x_1, x_2, \ldots, x_m$ of $\pi_1(D^2-Q_m)$. 
By $\mathcal{A}^b_*$, 
the braid group $B_m$ acts on $\pi_1(D^2-Q_m)$. 
It is presented by 
\[
\mathcal{A}^{\sigma_i}_*(x_j)=\begin{cases}
                        x_j x_{j+1} x_j^{-1} & \mathrm{if} \ j = i, \\
                        x_{j-1} & \mathrm{if} \ j = i+1, \\
                        x_j & \mathrm{otherwise}, 
\end{cases}
\]
and 
\[
\mathcal{A}^{\sigma_i^{-1}}_*(x_j)=\begin{cases}
                        x_{j+1} & \mathrm{if} \ j = i, \\
                        x_{j}^{-1} x_{j-1} x_{j} & \mathrm{if} \ j = i+1, \\
                        x_j & \mathrm{otherwise}, 
\end{cases}
\]
where $i=1,2,\ldots,m-1$ and $j=1,2,\ldots,m$. 

 \begin{prop} \label{Lem4-1}
For commutative $m$-braids $a$ and $b$, 
the link group of $\mathcal{S}_m(a,b)$ is presented by 
\begin{equation*} 
\pi_1(\mathbb{R}^4-\mathcal{S}_m(a,b))=\langle \, x_1 \,, \ldots,\, x_m \mid 
x_j=\mathcal{A}^a_*(x_j)=\mathcal{A}^b_* (x_j)\, ,\,\mathrm{for} \  j =1,2,\ldots,m \, \rangle. 
\end{equation*}
\end{prop}

\Proof 
Put $S=\mathcal{S}_m(a,b)$. 
Let $Q_m$ be a set of $m$ distinct interior points of $D^2$, 
and let $q_0$ be a point of $\partial D^2$. 
The space $N(T)-S$ 
is a fiber bundle over $T$ with the fiber $D^2-Q_m$ 
whose monodromy is given by 
 $\mathcal{A}^a$ and $\mathcal{A}^b$. 
Let us take commutative $\mathcal{A}^a$ and $\mathcal{A}^b$. 
Then we have  
  \begin{eqnarray*}
N(T)-S & \cong & ( (D^2-Q_m) \times I \times I ) / (x,0,u) \sim (\mathcal{A}^a(x),1,u), 
              (x,u,0) \sim (\mathcal{A}^b(x),u,1) \\
      &=& M \times I/(x,u,0) \sim (\mathcal{A}^b(x),u,1), 
\end{eqnarray*}
where 
\[
 M=((D^2-Q_m) \times I) / (x,0) \sim (\mathcal{A}^a(x),1), 
\]
and $x \in D^2-Q_m$ and $u \in I$. 

We compute $\pi_1(M)$, as follows. 
Since $M$ is a mapping torus whose monodromy is given by $\mathcal{A}^a$, by 
van Kampen's theorem, we can see that 
$\pi_1(M)$ has a presentation obtained from 
 $\pi_1(D^2-Q_m) * \mathbb{Z}$ by adding relations $s^{-1}x s= \mathcal{A}^a_*(x)$, 
where $x \in \pi_1(D^2-Q_m)$, and $s$ is the generator of $\mathbb{Z}$, 
which is represented by the loop $\{q_0\} \times S^1$ with the base point $q_0$. 
 Since $\pi_1(D^2-Q_m)$ is a free group 
generated by the standard generators $x_1, x_2, \ldots,x_m$ with the base point $q_0$, we have 
\begin{equation} \label{mm}
\pi_1(M)=\langle \, x_1, x_2, \dots,x_m, s \, \mid 
s^{-1} x_j s= \mathcal{A}^a_*(x_j) ,\,\mathrm{for} \  j =1,2,\ldots,m \rangle. 
\end{equation}

We compute $\pi_1(N(T)-S)$, as follows. 
Since 
$N(T)-S=M \times I/ (x,u,0) \sim (\mathcal{A}^b(x),u,1)$, 
where $x \in D^2-Q_m$ and $u \in I$ with 
$(x,u) \in M$, 
$N(T)-S$ 
is a mapping torus whose monodromy is given by $\mathcal{A}^b \times \mathrm{id}$. 
Thus we can see that 
$N(T)-S$ has a presentation obtained from 
 $\pi_1(M) * \mathbb{Z}$ by adding relations 
$t^{-1} y t=  (\mathcal{A}^b_* \times \mathrm{id}_*) (y)$, 
where $y \in \pi_1(M)$ and $t$ is the generator of $\mathbb{Z}$. 
Hence together with (\ref{mm}), we can see that 
$\pi_1(N(T)-S)$ is presented by 
\begin{eqnarray} \label{0611-6} 
&& \\
\lefteqn{\pi_1(N(T)-S)} && \nonumber \\
& =
\langle \, x_1, x_2, \dots,x_m, s, t \, \mid & 
s^{-1} x_j s= \mathcal{A}^a_*(x_j),\, 
t^{-1} x_j t= \mathcal{A}^b_*(x_j), \nonumber \\
 && {} t^{-1} s t= s,\,\mathrm{for} \  j =1,2,\ldots,m \,\rangle,  \nonumber 
\end{eqnarray}
where 
$s$ (resp. $t$) is represented by 
the loop $\mathbf{m}$ (resp. 
$\mathbf{l}$). 

We compute $\pi_1(S^4-S)$, as follows. 
We have $S^4-S=(N(T)-S) \cup _{\partial N(T)} E^4$. 
The fundamental group $\pi_1(N(T)-S)$ has the presentation (\ref{0611-6}). 
We obtain $\pi_1(E^4)$, as follows. 
Since $N(T)$ is a tubular neighborhood of $T$, and $T$ is 
the boundary of the standard unknotted solid torus in $\mathbb{R}^3 \times \{0\}$, we can see that 
the fundamental group of $E^4=\mathrm{cl}(S^4-N(T))$ is the knot group of a trivial torus knot. 
Hence $\pi_1(E^4)$ is an infinite cyclic group, 
where the generator $r$ is represented by the loop $\mathbf{r}$ 
(see \cite{Carter-Saito}, Section 5.2). 
Next we obtain $\pi_1(\partial E^4)$, as follows. 
 Since $\partial E^4=\partial N(T)=\partial D^2 \times T$ 
is a 3-dimensional torus $S^1 \times S^1\times S^1$, 
$\pi_1(\partial N(T))$ is 
isomorphic to $\mathbb{Z} \oplus \mathbb{Z} \oplus \mathbb{Z}$, 
where the generators $r^\prime$, 
$m^\prime$ and $l^\prime$ are 
 represented by the loops $\mathbf{r}$, $\mathbf{m}$ and $\mathbf{l}$ 
 respectively. 
Let $i_1 \,:\, \partial N(T) \rightarrow E^4$ and 
    $i_2 \,:\, \partial N(T) \rightarrow N(T)-S$ be 
inclusion maps. 
Since $i_{1*}(r^\prime)=r$, $i_{1*}(m^\prime)=1$, and $i_{1*}(l^\prime)=1$ 
in $\pi_1(E^4)$, and 
$i_{2*}(r^\prime)=x_1 x_2 \cdots x_m$, $i_{2*}(m^\prime)=s$, and $i_{2*}(l^\prime)=t$ 
 in 
 $\pi_1(N(T)-S)$, 
by van Kampen's theorem 
 $\pi_1(S^4-S)=\pi_1((N(T)-S) \cup _{\partial N(T)} E^4)$ is presented by 
\begin{eqnarray*}
& \langle \, r, x_1 \,, \ldots,\, x_m, s, t \mid &
 s^{-1} x_j s= \mathcal{A}^a_*(x_j),\, 
 t^{-1} x_j t = \mathcal{A}^b_*(x_j), t^{-1} s t = s, \\
&& {} r=x_1 x_2 \cdots x_m, 
s=1, t=1, 
\,\mathrm{for} \  j =1,2,\ldots,m \, \rangle, 
\end{eqnarray*}
which is the required formula. 
\qed
 \\

\begin{figure}
     \includegraphics*{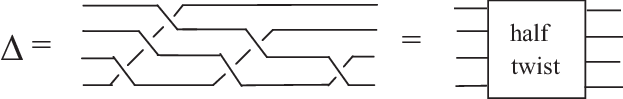}
  \caption{A half twist $\Delta$}
  \label{halftwist}
\end{figure}

   \begin{thm} \label{Thm4-6}
Put $S_n=\mathcal{S}_4(\sigma_1 \sigma_3, \Delta^{2n})$, 
where $\Delta= \sigma_1 \sigma_2 \sigma_3 \sigma_1 \sigma_2 \sigma_1$ 
(see Fig. \ref{halftwist}), 
and $n$ is a positive 
integer. 
Then the link group of $S_n$ is not a classical link group. 
Moreover, $S_n$ and $S_m$ are not equivalent for $n \neq m$. 
\end{thm}

Since $\Delta^{2n}$ is a central element of $B_4$ and the closure of $\sigma_1 \sigma_3$ 
is a classical link with two components, 
the torus-covering $T^2$-link $S_n$ has two components. 
Each component of $S_n$ is equivalent to $\mathcal{S}_2(\sigma_1, \sigma_1^{2n})$, 
which is equivalent to $\mathcal{S}_2(\sigma_1, e)$ by Corollary \ref{0629-2}. 
By Proposition \ref{eg}, it is equivalent to the spun $T^2$-knot of 
$\hat{\sigma_1}$. 
Since $\hat{\sigma_1}$ is a trivial knot, 
this is an unknotted $T^2$-knot. Thus 
each component of $S_n$ is an unknotted $T^2$-knot. 
\\

\Proof
By Proposition \ref{Lem4-1}, the link group $G_n$ of $S_n$ is computed as follows. 
Let $x_1\,, \ldots,\, x_4$ be the generators. 
The relations concerning the basis braid $\sigma_1 \sigma_3$ are 
 $x_1 = x_2$ and 
$x_3 = x_4$. 
The other relations concerning the basis braid $\Delta^{2n}$ are 
\begin{eqnarray*}
x_1 &=& (x_1 x_2 x_3 x_4)^{n}\, x_1\, (x_1 x_2 x_3 x_4)^{-n}, \\
x_2 &=& (x_1 x_2 x_3 x_4)^{n}\, x_2\, (x_1 x_2 x_3 x_4)^{-n}, \\
x_3 &=& (x_1 x_2 x_3 x_4)^{n}\, x_3\, (x_1 x_2 x_3 x_4)^{-n}, \\
x_4 &=& (x_1 x_2 x_3 x_4)^{n}\, x_4\, (x_1 x_2 x_3 x_4)^{-n} .
\end{eqnarray*}
Putting $a=x_1=x_2$ and $b=x_3=x_4$, we have  
 \[
G_n = \langle \, a\,,\, b \mid (a^2 b^2)^n b = b (a^2 b^2)^n \, ,\, 
(a^2 b^2)^n a=a (a^2 b^2)^n \, \rangle.
 \]
 
By Lemma \ref{0702-1}, $G_1$ is not a classical link group. 
Let us consider the case for $n>1$. 
For $n>1$, let $Z_n$ be the subgroup of $G_n$ generated by $h_n=(a^2 b^2)^n$. 
By Lemma \ref{center}, $Z_n$ is the center of $G_n$ for $n > 1$. 
Further, $Z_n$ is an infinite cyclic group; thus the center of $G_n$ is non-trivial. 

We will show that 
$G_n$ ($n>1$) is not a classical link group. 
Since the torus-covering $T^2$-link $S_n$ consists of two components, 
we show that 
$G_n$ is not a classical 2-component link group, as follows. 
 It is known \cite{Burde-Murasugi} that 
 if the center of the group of a classical 2-component link $L$ is non-trivial, 
then the link group of $L$ is isomorphic to one of the groups of type (a), (b), or (c) 
as follows: 
\begin{enumerate}[(a)]
\item $Z \times \mathbb{Z}$, 
\item $Z*_{Z} 
((Z \times \mathbb{Z})*_Z \mathbb{Z})$, 
\item $Z*_Z
\bigl( Z \times \mathbb{Z} )*_Z (\mathbb{Z}*_\mathbb{Z} \mathbb{Z}) \bigl)$,
\end{enumerate}
where 
 $\mathbb{Z}$ is an infinite cyclic group, and $Z=\langle \, h\, \rangle$ is 
 a \lq\lq special" infinite cyclic 
group which is the center of the link group of types (b) and (c). 
In type (b), the amalgamation concerning the last factor $\mathbb{Z}=\langle \, q\, \rangle$ is given by 
$h=q^\alpha$ for an integer $\alpha>1$. 
In type (c), the last factor $\mathbb{Z}*_\mathbb{Z} \mathbb{Z}$ is the group of the torus $(\alpha,\, \beta)$-knot, i.e. 
$\mathbb{Z}*_\mathbb{Z} \mathbb{Z}=\langle \, x\,,\, y\mid x^\alpha=y^\beta \, \rangle$ for 
coprime positive integers $\alpha$ and $\beta$, 
and the amalgamation is given by $h=x^\alpha=y^\beta$. 

Since an infinite cyclic group $Z_n$ is the center of $G_n$ by Lemma \ref{center}, 
it suffices to show that $G_n$ ($n>1$) is neither of type (b) nor (c). 
Further, in these cases, the center $Z$ equals $Z_n$ and $h=h_n$.  

{\bf (Case (b))} 
If $G_n$ is of type (b), then 
$G_n=(Z_n \times \mathbb{Z})*_{Z_n} \mathbb{Z}=
(\langle \,h_n\, \rangle \times \langle\, k\, \rangle)*_{Z_n} \langle \,q \, \rangle$, 
where the amalgamation is given by $h_n=q^\alpha$ for an integer $\alpha>1$. 
Put $h_n^\prime=f(h_n)$ and $q^\prime=f(q)$ for 
a natural epimorphism 

\begin{equation} \label{f}
f\,: \, G_n \rightarrow \mathbb{Z}/2 \mathbb{Z} * \mathbb{Z}/2 \mathbb{Z} =\langle \, a^\prime \, \rangle * 
\langle \, b^\prime \, \rangle, 
\end{equation}
where 
$a^\prime=f(a)$ and $b^\prime=f(b)$, which are the basis. 
Since $h_n=(a^2b^2)^n$, we see that $h_n^\prime=1$. 
Since  $h_n=q^\alpha$, it follows that $q^{\prime \, \alpha}=1$. 
If $q^\prime=1$, then $f(G_n)=\langle \,f(k) \, \rangle$, which 
is generated by at most one generator. 
However, $f(G_n)$ is generated by two generators $a^\prime$ and $b^\prime$. 
Hence $q^\prime$ is non-trivial. 
Since $\alpha>1$ and a non-trivial element of $f(G_n)=
\mathbb{Z}/2 \mathbb{Z} * \mathbb{Z}/2 \mathbb{Z}$ 
has order $2$ or $\infty$ by Lemma \ref{lem1}, $q^\prime$ has order $2$; thus $\alpha=2$. 
Consider the abelianization map 
\begin{equation} \label{0824-2}
\phi\, :\, G_n \rightarrow G_n/[G_n, G_n]=\mathbb{Z} \times \mathbb{Z}
\end{equation}
 and put 
 $\overline{a}=\phi(a)$ and $\overline{b}=\phi(b)$, which are the basis. 
 Since $h_n=q^2$ and $\phi(h_n)=\overline{a}^{2n} \overline{b}^{2n}$, 
 it follows that $\phi(q)=\overline{a}^{n} \overline{b}^{n}$. 
Consider the abelianization map 
\begin{equation*} \label{0824-1}
\phi^\prime \, :\, f(G_n)=\mathbb{Z}/2 \mathbb{Z} * \mathbb{Z}/2 \mathbb{Z} 
 \rightarrow \mathbb{Z}/2 \mathbb{Z} 
\times \mathbb{Z}/2 \mathbb{Z}, 
\end{equation*}
and put $\overline{a}^\prime=\phi^\prime(a^\prime)$ and 
$\overline{b}^\prime=\phi^\prime(b^\prime)$, which are the basis. 
Since $\phi(q)=\overline{a}^{n} \overline{b}^{n}$, it follows that  
$\phi^\prime(q^\prime)=1$ if $n$ is even, and 
 $\phi^\prime(q^\prime)=\overline{a}^\prime \overline{b}^\prime$ if $n$ is odd; thus 
$\phi^\prime(q^\prime)=1$ or $\overline{a}^\prime \overline{b}^\prime$. 
However, 
$\phi^\prime(q^\prime)=\overline{a}^\prime$ or $\overline{b}^\prime$, as follows. 
Since $h_n=q^2$, it follows that $q^{\prime \, 2}=1$; thus 
$q^\prime=\xi^{-1} a^\prime \xi$ or $\xi^{-1} b^\prime \xi$ for some 
$\xi \in f(G_n)$ by 
Lemma \ref{lem1}. Thus 
$\phi^\prime(q^\prime)=\overline{a}^\prime$ or $\overline{b}^\prime$. 
This is 
 a contradiction. 

{\bf (Case (c))} 
If $G_n$ is of type (c), then 
$G_n=(Z_n \times \mathbb{Z})*_{Z_n} (\mathbb{Z}*_\mathbb{Z} \mathbb{Z})$, 
where $\mathbb{Z}*_\mathbb{Z} \mathbb{Z}=\langle \, x\,,\, y\mid x^\alpha=y^\beta \, \rangle$ for 
coprime positive integers $\alpha$ and $\beta$, 
and the amalgamation is given by $h_n=x^\alpha=y^\beta$. 
Since $h_n=(a^2b^2)^n$, we see that $h_n^\prime=1$, 
where $h_n^\prime=f(h_n)$. 
Since  $h_n=x^\alpha=y^\beta$, it follows that $x^{\prime \, \alpha}=y^{\prime \, \beta}
= 1$, where $x^\prime=f(x)$ and $y^\prime=f(y)$. 
If $x^\prime=1$ and $y^\prime=1$, then $f(\langle \, x\,,\, y \mid x^\alpha=y^\beta \, \rangle)=1$ 
and it follows that $f(G_n)$ 
is generated by at most one generator. 
However, $f(G_n)$ is generated by two generators $a^\prime$ and $b^\prime$. 
Hence we can assume that $x^\prime$ is non-trivial. 
Since any element of $f(G_n)=\mathbb{Z}/2 \mathbb{Z} * \mathbb{Z}/2 \mathbb{Z}$ 
has order $2$ or $\infty$ by Lemma \ref{lem1}, 
it follows that 
  $\alpha=2$. 
Since $h_n=x^2$ and $\phi(h_n)=\overline{a}^{2n} \overline{b}^{2n}$, 
it follows that $\phi(x)=\overline{a}^{n} \overline{b}^{n}$, and hence $\phi^\prime(x^\prime)=1$ or 
 $\overline{a}^\prime \overline{b}^\prime$ by the same argument as in Case (b). 
However, since $h_n=x^2$, it follows that $x^{\prime \, 2}=1$, and hence 
$\phi^\prime(x^\prime)=\overline{a}^\prime$ or $\overline{b}^\prime$ 
  by the same argument as in Case (b). 
 This is a contradiction. 
Thus $G_n$ is not a classical link group. 
 
Now we will show that 
$S_n$ and $S_m$ are not equivalent for $n \neq m$. The center $Z_1$ of $G_1$ 
is a free abelian group of rank $2$ (see the proof of Lemma \ref{0702-1}), while, for $n>1$, 
the center $Z_n$ of $G_n$ is an infinite cyclic group by Lemma \ref{center}; 
thus 
it suffices to show in the case when $n, m>1$. 
The abelianization of $G_n/Z_n$ for $n>1$ is $\mathbb{Z} \times \mathbb{Z}/2n \mathbb{Z}$;
thus $G_n$ is not isomorphic to $G_m$ for $n \neq m$, and hence 
$S_n \not\sim S_m$ for $n \neq m$.  
%
%
\qed
\\

A {\it 2-link} is a surface link whose components are homeomorphic to 2-spheres. 
It is known \cite[Chapter 3, Corollary 2]{Hillman} that if the center of a $\mu$-component 2-link group with $\mu>1$ is non-trivial, then 
the center must be a torsion group. Hence we have a corollary. 

 \begin{cor} \label{Cor4-6}
 The link group of the 2-component torus-covering $T^2$-link of Theorem \ref{Thm4-6} 
 is not a 2-component 2-link group. 
 \end{cor}
 
 \Proof 
  For any $n>0$, the center $Z_n$ of $G_n$ is 
 non-trivial and torsion free by Lemmas \ref{0702-1} and \ref{center}; 
thus the conclusion follows from \cite[Chapter 3, Corollary 2]{Hillman}. 
\qed 
 \\

\begin{lem} \label{0702-1}
The group $G_1$ of Theorem \ref{Thm4-6} is not a classical link group.  
\end{lem}
\Proof
Let $Z_1$ be the subgroup of $G_1$ generated by $\{ \, a^2 \,,\, b^2 \, \}$. 
We will show that $Z_1$ is the center of $G_1$, as follows. 
Let $N$ be a normal subgroup of $G_1$. If the center of $G_1/N$ is trivial, 
then $N$ contains the center of $G_1$. 
Since $Z_1$ consists of central elements, $Z_1$ is a normal subgroup of $G_1$ 
such that $Z_1$ is contained in the center of $G_1$. 
Hence 
it suffices to show that the center of the quotient group $G_1 /  Z_1$ is trivial. 
Since 
$G_1 / Z_1 = \mathbb{Z}/2 \mathbb{Z} * \mathbb{Z}/2 \mathbb{Z}$, 
the center of $G_1/Z_1$ is trivial; thus $Z_1$ is the center of $G_1$. 
Let us take the abelianization map $\phi$ 
given by (\ref{0824-2}). 
 Since $\phi(a^2)=\overline{a}^2$ and $\phi(b^2)=\overline{b}^2$, the center $Z_1$ is a 
 free abelian group of rank 2; thus 
$Z_1$ is generated by 
two generators. Hence it follows from \cite{Burde-Murasugi} that if $G_1$ is a classical link group, then it is 
isomorphic to $\mathbb{Z} \times \mathbb{Z}$ (type (a)): $G_1$ is commutative. 
However, since the image of $G_1$ by the natural epimorphism $f$ given by (\ref{f}) 
is a non-commutative group $\mathbb{Z}/2 \mathbb{Z} * \mathbb{Z}/2 \mathbb{Z}$, 
$G_1$ is not commutative. This is a contradiction. 
Thus $G_1$ is not a classical link group.  
\qed

\begin{lem} \label{center}
Let us consider the group $G_n$ of Theorem \ref{Thm4-6}. 
For $n>1$, let $Z_n$ be the subgroup of $G_n$ 
generated by $h_n=(a^2 b^2)^n$. 
Then $Z_n$ is the center of $G_n$. 
Moreover $Z_n$ is an infinite cyclic group. 
\end{lem}
\Proof 
By the same argument as in the proof of Lemma \ref{0702-1}, 
 in order to show that $Z_n$ is the center of $G_n$, 
it suffices to show that the center of the quotient group $G_n /  Z_n$ is trivial. 
We see that  
$G_n / Z_n = \langle \, a \,,\, b \mid (a^2 b^2)^n=1 \, \rangle$, which is an amalgamated product 
$\langle \, a \, \rangle *_U \langle \, b \,,\, x \mid x^n=1 \, \rangle$, 
where $U=\langle \, a^2 \, \rangle=\langle \, x b^{-2} \, \rangle=\mathbb{Z}$ and 
the amalgamation is given by $a^2 =x b^{-2}$. 
Put $H_1=\langle \,a\, \rangle$ and $H_2=\langle \, b \,,\, x \mid x^n=1 \, \rangle$. 
We can take $\{ \,1\,,\,a \, \}$ 
as a set of right-handed coset representatives of $U$ in $H_1$. 

Let $h$ be a central element of $G_n / Z_n=H_1*_U H_2$. 
By \cite{Neumann} or  
\cite[p.73, Theorem 11.3]{Bogopoloski}, 
$h$ is uniquly written as  
$h=u a^{\delta} c_1 a c_2 \cdots a c_t a^{\epsilon}$, 
where $u \in U$ and $c_1\,, \ldots,\, c_t$ are non-trivial elements of 
a set of right-handed coset representatives of $U$ 
in $H_2$, and 
$\delta, \epsilon \in \{0,1 \}$, 
which is called a {\it normal form}. 
 Since $ah=ha$, it follows that  
 $a u a^{\delta} c_1 a c_2 \cdots a c_t a^{\epsilon} =
 u a^{\delta} c_1 a c_2 \cdots a c_t a^{\epsilon} a$. 
 Since $ua=au$ in the amalgamated product $H_1*_U H_2$, it follows that 
 $ u a^{\delta} a c_1 a c_2 \cdots a c_t a^{\epsilon} =
 u a^{\delta} c_1 a c_2 \cdots a c_t a a^{\epsilon}$, and hence 
 $a c_1 a c_2 \cdots a c_t  =
  c_1 a c_2 \cdots a c_t a$ as elements in $H_1*_U H_2$. 

 If $t >0$, then $ac_1 a c_2 \cdots a c_t$ and $c_1 a c_2 \cdots a c_t a$ are in 
 distinct normal forms, which is a contradiction. 
 Hence $t=0$ and 
 $h=u a^{\delta}=a^k$ 
  for an integer $k$. 
Since $hb=bh$, $a^k b=b a^k$. 
If $k=1$, then $ab=ba$. 
In this case, if $b$ is not in $U$, then we can 
take $b$ as a non-trivial right-handed coset representative of $U$ in $H_2$. 
 It follows that then $ab$ and $ba$ are in distinct normal forms, 
which is a contradiction. 
If $k=2l+1$ (resp. $k=2l$) for a non-zero integer $l$, then 
    $a^k b=uab$ and $b a^k=c a$ (resp. $a^k b=ub$ and $b a^k= c$), 
where in both cases $u=a^{2l} \in U$ and $c=b (xb^{-2})^l$. 
In these cases, if neither $b$ nor $c$ is in $U$ and we can take $b$ and $c$ 
as distinct right-handed coset 
representatives of $U$ in $H_2$, 
then $a^k b$ and $b a^k$ have distinct normal forms $uab$ and $ca$ (resp. $ub$ and $c$), 
which is a contradiction. 
Then it follows that $k=0$ and hence $h=1$; thus the center of $G_n /  Z_n$ is trivial. 

It remains to show that for a non-zero 
integer $l$, neither $b$ nor $c=b (xb^{-2})^l$ is in $U$ and 
we can take $b$ and $c$ as distinct right-handed coset 
representatives of $U$ in $H_2$. 
%
The group $H_2=\langle \, b \,,\, x \mid x^n=1 \, \rangle$ is the free product 
of $\langle \,b\, \rangle$ and $\langle \,x \mid x^n=1 \, \rangle$. 
By \cite{Neumann} or  
\cite[p.73, Theorem 11.3]{Bogopoloski} again, every element of $H_2$ has a normal 
form $b_1^{\delta} x_1 b_2 x_2 \cdots b_t x_t^{\epsilon}$, 
where $b_1, b_2,\ldots, b_t$ (resp. $x_1,x_2, \ldots, x_t$) are non-trivial elements 
of $\langle b \rangle$ (resp. $\langle \,x \mid x^n=1 \, \rangle$), 
and $\delta, \epsilon \in \{0,1\}$. 
Let us determine the normal forms of $b$ and $c$. 
Put $l_0=|l|$, a positive integer. 
We can see that $b$ has a normal form $b$, and $c=b (xb^{-2})^l$ has a normal form 
$b(xb^{-2})^{l_0}$ (resp. $b^3 x^{-1} (b^2x^{-1})^{l_0-1}$) if $l>0$ (resp. $l<0$). 
Further, an element of $U=\langle \, xb^{-2} \, \rangle$ in $H_2$ 
has a normal form $1$, 
$(xb^{-2})^{m_0}$ or $(b^2x^{-1})^{m_0}$, where $m_0$ is a positive integer. 
Hence, by the uniqueness of normal forms, we can see that neither $b$ nor $c$ is in $U$. 
Similarly, 
 an element of $Ub$ has a normal form $b$, $(x b^{-2})^{m_0-1} x b^{-1}$ or 
$(b^2 x^{-1})^{m_0} b$, and 
an element of $Uc=Ub (xb^{-2})^l$ has a normal form 
$b(xb^{-2})^{l_0}$, 
$b^3 x^{-1} (b^2x^{-1})^{l_0-1}$, 
$(x b^{-2})^{m_0-1} x b^{-1} (xb^{-2})^{l_0}$, 
$(x b^{-2})^{m_0-1} x b x^{-1} (b^2 x^{-1})^{l_0-1}$, 
$(b^2 x^{-1})^{m_0} b (xb^{-2})^{l_0}$, or 
$(b^2 x^{-1})^{m_0} b^3 x^{-1} (b^2 x^{-1})^{l_0-1}$. 
By the uniqueness of normal forms, we can see that $Ub \neq Uc$. 
Thus, for a non-zero 
integer $l$, neither $b$ nor $c$ is in $U$ and 
we can take $b$ and $c$ as distinct right-handed coset 
representatives of $U$ in $H_2$, 
and it follows that the center of $G_n /  Z_n$ is trivial. 
Therefore $Z_n$ is the center of $G_n$. 

 Let us take the abelianization map $\phi$ given by (\ref{0824-2}).  
  Since the image $\phi(h_n)$ is $\overline{a}^{2n} \overline{b}^{2n}$,  
 the center $Z_n$ ($n>1$) is an infinite cyclic group. 
\qed

\begin{lem} \label{lem1}
Let us consider the group 
$\mathbb{Z}/p \mathbb{Z} *\mathbb{Z}/q \mathbb{Z}$, the free product of 
$\mathbb{Z}/p \mathbb{Z}$ and 
$\mathbb{Z}/q \mathbb{Z}$, where 
$p$, $q$ are positive integers greater than one. 
Then the order of a non-trivial element 
$\mathbb{Z}/p \mathbb{Z} *\mathbb{Z}/q \mathbb{Z}$ is either $\infty$, 
a divisor of $p$, or a divisor  
of $q$. 
Further, if the order is finite, then 
the element can be written as a conjugate of 
an element of the same order in $\mathbb{Z}/p \mathbb{Z}$ or 
$\mathbb{Z}/q \mathbb{Z}$. 
\end{lem}

\Proof
Let $z$ be a non-trivial element of $\mathbb{Z}/p \mathbb{Z} * \mathbb{Z}/q \mathbb{Z}$. 
By \cite{Neumann} or  
\cite[p.73, Theorem 11.3]{Bogopoloski}, $z$ 
has a normal form 
\begin{eqnarray} \label{normal-z}
&&  x_1 y_1 x_2 y_2 \cdots x_t y_t, \label{1-1} \\
&&  x_1 y_1 x_2 y_2 \cdots x_{t-1} y_{t-1} x_t,  \label{2} \\
&&  y_1 x_2 y_2 \cdots x_t y_t, \label{3} \\
&&  y_1 x_2 y_2 \cdots x_{t-1} y_{t-1} x_t, \label{4}
\end{eqnarray}
where $t$ is an integer with $t>1$ in (\ref{4}) and $t>0$ otherwise, and 
$x_1, x_2,\ldots, x_t$ (resp. $y_1, y_2,\ldots, y_t$) 
are non-trivial elements of $\mathbb{Z}/p \mathbb{Z}$ (resp. $\mathbb{Z}/q \mathbb{Z}$). 
 
In cases (\ref{1-1}) and (\ref{4}), 
 $z^l$ has a normal form which is not $1$ 
 for any positive integer $l$. 
Thus the order of $z$ is infinite. 

  In cases (\ref{2}) and (\ref{3}), we show the lemma 
using induction for $t$, as follows. 
If $t=1$, then 
$z=x_1$ (resp. $y_1$) for (\ref{2}) (resp. (\ref{3})) and 
the order of $z$ is a divisor of $p$ (resp. $q$). 
Now let us assume that if $t<s$, then the order of any element $z$ with the normal form 
(\ref{2}) or (\ref{3}) 
is infinite or a divisor of $p$ or $q$. 
Let us consider $z$ with the normal form (\ref{2}) with $t=s$. 
If $x_1 x_t = 1$, then $x_1^{-1} z x_1$ has a normal form (\ref{3}) with $t=s-1$; 
thus, by the assumption, the statement of the lemma holds. 
If $x_1 x_t \neq 1$, then $z^l$ has a non-trivial normal form 
 for any positive integer $l$; thus the order of $z$ is infinite. 
For $z$ with the normal form (\ref{3}), we can apply the same argument. 
Further we can see that if the order $\mathrm{ord}(z)$ 
of $z$ is finite, i.e. a divisor of $p$ or $q$, then 
 $z$ can be written as 
$\xi^{-1} x \xi$, where 
$\xi \in \mathbb{Z}/p \mathbb{Z} * \mathbb{Z}/q \mathbb{Z}$, and 
 $x$ is an element of order $\mathrm{ord}(z)$ in 
$\mathbb{Z}/p \mathbb{Z}$ or $\mathbb{Z}/q \mathbb{Z}$. 
\qed
\\

We can consider the knot version of Theorem \ref{Thm4-6}. 
 
\begin{thm} \label{Thm4-9}
Put $S_n=\mathcal{S}_4(\sigma_1 \sigma_3, \Delta^{2n+1})$, 
where $n$ is a positive integer. 
Then the knot group of $S_n$ is not a classical knot group. 
 Moreover, $S_n$ and $S_m$ are not equivalent for $n \neq m$. 
\end{thm}

The torus-covering $T^2$-knot 
$S_0=\mathcal{S}_4( \sigma_1 \sigma_3, \Delta)$ is unknotted 
(see Corollary \ref{0704-6}). 
\\

\Proof 
By Proposition \ref{Lem4-1}, the knot group $G_n$ of $S_n$ is computed as follows. 
Let $x_1, \ldots, x_4$ be the generators. 
Then the relations concerning the basis braid $\sigma_1 \sigma_3$ are 
$x_1 = x_2$ and 
$x_3 = x_4$. 
The other relations concerning the basis braid $\Delta^{2n+1}$ are 
\begin{eqnarray*}
x_1 &=& (x_1 x_2 x_3 x_4)^{n}\, x_1 x_2 x_3 x_4 x_3^{-1} x_2^{-1} x_1^{-1}\, (x_1 x_2 x_3 x_4)^{-n}, \\
x_2 &=& (x_1 x_2 x_3 x_4)^{n}\, x_1 x_2 x_3 x_2^{-1} x_1^{-1}\, (x_1 x_2 x_3 x_4)^{-n}, \\
x_3 &=& (x_1 x_2 x_3 x_4)^{n}\,  x_1 x_2 x_1^{-1}\, (x_1 x_2 x_3 x_4)^{-n}, \\
x_4 &=& (x_1 x_2 x_3 x_4)^{n}\,  x_1\, (x_1 x_2 x_3 x_4)^{-n}.
\end{eqnarray*}
Putting $a=x_1=x_2$ and $b=x_3=x_4$, we have  
 \[
G_n = \langle \, a\,,\, b \mid b(a^2 b^2)^n =(a^2 b^2)^n a\,, \, a(a^2 b^2)^{n+1} = (a^2 b^2)^{n+1} b \, \rangle .
\]

Let us assume that $G_n$ is a classical knot group. 
Let $Z_n$ be the subgroup of $G_n$ generated by $h_n=(a^2 b^2)^{2n+1}$, which is a central element. 
By Lemma \ref{center2}, $Z_n$ is the center of $G_n$. 
Further, $Z_n$ is an infinite cyclic group. 
It is known \cite{Burde-Zieschang} that if the center of a classical knot group is non-trivial, 
then the knot is a torus knot. 
Hence, by the assumption, 
$G_n$ is isomorphic to a torus knot group. 
Let $G_{p,\,q}$ be the $(p,\,q)$-torus knot group isomorphic to $G_n$, where 
$p$ and $q$ are coprime positive integers. 
Let $Z_{p,\,q}$ be the center of $G_{p,\,q}$. 
Then $G_{p,\,q}=\langle \,x,\,y \mid x^p=y^q \, \rangle$ and 
$Z_{p,\,q}$ is generated by $h=x^p=y^q$. 
Put $G_{p,\,q}^\prime=G_{p,\,q}/Z_{p,\,q}$, which is 
$\langle \,x,\,y \mid x^p=y^q=1\, \rangle= \mathbb{Z}/p \mathbb{Z} * \mathbb{Z}/q \mathbb{Z}$. 
The abelianization of $G_{p,\,q}^\prime$ 
is isomorphic to 
$\mathbb{Z}/p \mathbb{Z} \times \mathbb{Z}/q \mathbb{Z}$. 
%
   %

Consider the quotient group 
$G_n^\prime=G_n/Z_n$. 
By (\ref{0704-4}) in the proof of Lemma \ref{center2}, 
$G_n^\prime=\langle \, a\,,\, x \mid x^{2n+1}=(a^2 x^n)^2=1 \, \rangle$. 
   The abelianization of $G_n^\prime$ is presented by 
   $\langle \, a\,,\, x \mid x^{2n+1}=(a^2 x^n)^2=1, \, ax=xa \, \rangle$, which equals 
   $\langle \, a \mid a^{4(2n+1)} \, \rangle=
   \mathbb{Z}/4(2n+1)\mathbb{Z}$. 
Since $G_{p,\,q}$ and $G_n$ are isomorphic, 
so are the abelianizations of $G_{p,\,q}^\prime$ and 
$G_n^\prime$. 
Hence, comparing the order of the groups we see that $p\,q=4(2n+1)$. 
Since $G_n^\prime$ has an element of order $2n+1$ by Lemma \ref{0704-5}, 
and the order of a non-trivial torsion element of 
$G_{p,\,q}^\prime$ is a divisor of $p$ or $q$ by Lemma \ref{lem1}, 
it follows that $2n+1$ is a divisor of $p$ or $q$. Hence we can determine coprime positive integers $p$ and $q$ 
by $p=4$ and $q=2n+1$. 

For any element $z$ of order $2$ in 
$G_{4,\,2n+1}^\prime=\mathbb{Z}/4 \mathbb{Z} *\mathbb{Z}/(2n+1) \mathbb{Z}$, 
$z$ can be written as $z=z^{\prime \, 2}$ for some element $z^\prime$ of order $4$ by 
Lemma \ref{0704-2}.  
   %
Since $y=a^2 x^{n} \in G_n^\prime$ is of order $2$ by Lemma \ref{0704-5}, 
and $G_n^\prime$ and $G_{4,\,2n+1}^\prime$ 
are isomorphic, there is an element $y^\prime \in G_n^\prime$ with $y=y^{\prime \, 2}$, 
and hence  
$G_n^\prime
= \langle \, a\,,\, x\,,\, y^\prime \mid x^{2n+1}=1\,,\, y^{\prime \, 4}=1\,,\, a^2 x^n=y^{\prime \, 2} \, \rangle$. 
Let $N_w$ be the normal subgroup of $G^\prime_{4,\,2n+1}$ generated by an element $w$ 
of order $2n+1$. 
The quotient group $G_{4,\,2n+1}^\prime/N_w$ does not depend on the choice 
of $w$ and 
$G_{4,\,2n+1}^\prime/N_w=\mathbb{Z}/4 \mathbb{Z}$ 
by Lemma \ref{0704-3}. 
%
 We will denote it by 
$G_{4,\,2n+1}^\prime/N$. 
Let $N_x$ be the normal subgroup of $G_n^\prime$ generated by $x$. 
Since $x$ has order $2n+1$ by Lemma \ref{0704-5}, 
$G_n^\prime/N_x$ is isomorphic to 
$G_{4,\,2n+1}^\prime/N=\mathbb{Z}/4 \mathbb{Z}$: $G_n^\prime/N_x$ is abelian. 
Adding the relation $x=1$ to the presentation of $G_n^\prime$, we see that 
 $G_n^\prime/N_x 
= \langle \, a\,,\, y^\prime \mid y^{\prime \, 4}=1\,,\, a^2 =y^{\prime \ 2} \, \rangle$. 
 Since there is a natural epimorphism $f$ from $G_n^\prime/N_x$ onto 
 $\mathbb{Z}/2 \mathbb{Z} * \mathbb{Z}/2 \mathbb{Z}$ with the basis $f(a)$ and $f(y^\prime)$, 
 $G_n^\prime/N_x$ is not abelian. 
This is a contradiction. 
Thus $G_n$ is not a classical knot group. 
Since the abelianization of $G_n^\prime$ is $\mathbb{Z}/4(2n+1) \mathbb{Z}$, 
it follows that $S_n \not\sim S_m$ for $n \neq m$. 
\qed

\begin{lem} \label{center2}
Let us consider the group $G_n$ of Theorem \ref{Thm4-9}. 
For $n>0$, 
let $Z_n$ be the subgroup of $G_n$ generated by $h_n=(a^2 b^2)^{2n+1}$. 
Then $Z_n$ is the center of $G_n$. Moreover 
$Z_n$ is an infinite cyclic group. 
\end{lem}
\Proof
Since $h_n$ is a central element, by the same argument as in the proof of Lemma \ref{0702-1}, 
in order to show that $Z_n$ is the center of $G_n$, 
it suffices to show that the center of $G_n^\prime=G_n/Z_n$ is trivial. 
The quotient group 
$G_n^\prime=G_n/Z_n$ is presented by 
$\langle \, a\,,\, b\,,\, x \mid x=a^2 b^2\,,\, b x^n =x^n a\,,\, a x^{n+1}=x^{n+1} b\,,\, x^{2n+1}=1 \, \rangle$. 
By eliminating $b$ by $b=x^n a x^{-n}$, we have   
\begin{equation} \label{0704-4}
G_n^\prime = \langle \, a\,,\, x \mid x^{2n+1}=(a^2 x^n)^2=1 \, \rangle, 
\end{equation}
which is an amalgamated product 
$\langle \, a \, \rangle *_U \langle \, x\,,\, y \mid x^{2n+1}=1\,, \, y^2=1 \, \rangle$, 
where $U=\langle \, a^2 \, \rangle=\langle \, yx^{-n} \, \rangle=\mathbb{Z}$ and 
the amalgamation is given by $a^2 =yx^{-n}$. 
We can show that the center of $G_n^\prime$ is trivial by 
the following argument similar to the proof of Lemma \ref{center}, 
as follows. 
Put $H_1=\langle \,a \,\rangle$ and 
$H_2=\langle \, x\,,\, y \mid x^{2n+1}=1 \,, \, y^2=1 \, \rangle$. 
Note that $H_2$ 
is a free product of $\langle \, x \mid x^{2n+1}=1 \, \rangle=\mathbb{Z}/(2n+1) \mathbb{Z}$ and 
$\langle \, y \mid y^2=1 \, \rangle=\mathbb{Z}/2 \mathbb{Z}$. 
By \cite{Neumann} or  
\cite[p.73, Theorem 11.3]{Bogopoloski}, 
any element of $G_n^\prime=H_1*_U H_2$ has a normal form 
$u a^{\delta} c_1 a c_2 \cdots a c_t a^{\epsilon}$, where 
$u \in U$ and $c_1, \ldots, c_t$ are non-trivial elements of 
a set of right handed coset representatives 
of $U$ in $H_2$, and $\delta, \epsilon \in \{0, 1\}$. 

Let $h$ be a central element of $G_n^\prime$. 
By the same argument as in the proof of Lemma \ref{center}, 
by using normal forms, 
we see that 
 $h=a^k$ 
  for an integer $k$. 
Since $hx=xh$, $a^k x=x a^k$. 
If $k=2l+1$ (resp. 2l) for a non-zero integer $l$, then  
 $a^k x =uax$ and 
 $x a^k = c a$ 
 (resp. $a^k x =ux$ and $x a^k = c$), 
where $u=a^{2l} \in U$ and $c= x(yx^{-n})^l$ in both cases. 
If neither $x$ nor $c$ is in $U$ and 
we can choose $x$ and $c$ as distinct right-handed coset representatives of $U$ in $H_2$, then 
in both cases $a^k x$ and $x a^k$ have distinct normal forms, 
which is a contradiction. Then $k=0$, and it follows that $h=1$. 

It remains to show that 
neither $x$ nor $c=x(yx^{-n})^l$ ($l \neq 0$) is in $U$ and 
we can choose $x$ and $c$ as distinct right-handed coset representatives of $U$ in $H_2$. 
By \cite{Neumann} or  
\cite[p.73, Theorem 11.3]{Bogopoloski}, 
any element of the free product 
$H_2=\langle \, x \mid x^{2n+1}=1 \, \rangle * \langle \, y \mid y^2=1 \, \rangle$ 
has a normal form 
$x_1^{\delta m_1} y x_2 y \cdots x_t y^\epsilon$, 
where $x_1, x_2,\ldots, x_t$ are non-trivial elements of 
 $\langle x \mid x^{2n+1}=1 \rangle$ and $\delta, \epsilon \in \{0,1\}$. 
Let us determine the normal forms of $x$ and $c$. 
Put $l_0=|l|$, a positive integer. 
We can see that $x$ has a normal form $x$, and  
$c=x(yx^{-n})^l$ has a normal form $x(yx^{-n})^{l_0}$ (resp. $x^{n+1}y(x^n y)^{l_0-1}$) 
if $l>0$ (resp. $l<0$). 
Further, 
an element of $U=\langle \, yx^{-n} \, \rangle$ has a normal form either $1$, 
$(yx^{-n})^{m_0}$ or $(x^n y)^{m_0}$, where $m_0$ is a positive integer. 
Hence, by the uniqueness of normal forms, neither $x$ nor $c$ is in $U$. 
Similarly, 
if $n=1$ (resp. $n>1$), then 
an element of $Ux$ has a normal form either $x$, 
$(yx^{-1})^{m_0-1} y$, or $(x y)^{m_0}x$ 
(resp.  $x$, 
$(yx^{-n})^{m_0-1} yx^{-n+1}$, or $(x^n y)^{m_0}x$). 
Hence in both cases $c$ is not an element of $Ux$. 
Thus neither $x$ nor $c$ is in $U$ and 
we can choose $x$ and $c$ as distinct right-handed coset representatives of $U$ in $H_2$, 
and it follows that the center of $G_n^\prime$ is trivial. 
Thus $Z_n$ is the center. 
Considering the abelianization map of $G_n$, we see that $Z_n$ is an infinite cyclic group. 
\qed

\begin{lem} \label{0704-5}
The element $x$ (resp. $y=a^2 x^{n}$) 
of $G_n^\prime$ of Theorem \ref{Thm4-9} (see (\ref{0704-4})) is of order $2n+1$ (resp. $2$). 
\end{lem}
\Proof
By Lemma \ref{center2}, $G_n^\prime$ is an amalgamated product. 
Seeing the normal forms of the powers of $x$ and $y$, 
 we can show that the order of $x$ is $2n+1$ and 
the order of $y$ is $2$. 
\qed

\begin{lem} \label{0704-2}
For any element $z$ of order $2$ in 
$\mathbb{Z}/4 \mathbb{Z} *\mathbb{Z}/(2n+1) \mathbb{Z}$, 
$z$ can be written as $z=z^{\prime \, 2}$ for some element $z^\prime$ of order $4$.  
\end{lem}
 \Proof
The order of $z$ is $2$. It is a divisor of $4$, 
and it is not a divisor of $2n+1$. Thus it follows from 
Lemma \ref{lem1} that 
$z=\xi^{-1} u^2 \xi$, 
where 
$\xi \in \mathbb{Z}/4 \mathbb{Z} * \mathbb{Z}/(2n+1) \mathbb{Z}$ and 
$u$ is a generator of $\mathbb{Z}/4 \mathbb{Z}$. 
Thus $z$ can be written as  
$z=z^{\prime \,2}$, where $z^\prime=\xi^{-1} u \xi$. 
Since the order of $u$ is 4, so is the order of $z^\prime$. 
\qed

\begin{lem} \label{0704-3}
Put $G_{4,\,2n+1}^\prime=\mathbb{Z}/4\mathbb{Z} * \mathbb{Z}/(2n+1)\mathbb{Z}$, 
and let $N_w$ be a normal subgroup of $G_{4,\,2n+1}^\prime$ generated by an 
element $w$ of order $2n+1$. 
Then $G_{4,\,2n+1}^\prime/N_w$ does not depend on the choice 
of $w$, and $G_{4,\,2n+1}^\prime/N_w=\mathbb{Z}/4 \mathbb{Z}$. 
\end{lem}
\Proof
Since $w$ has order $2n+1$, 
$w=\xi^{-1} 
 v^k
\xi$ by Lemma \ref{lem1},
where $\xi \in \mathbb{Z}/4 \mathbb{Z} * \mathbb{Z}/(2n+1) \mathbb{Z}$ and 
$v$ is a generator of $\mathbb{Z}/(2n+1) \mathbb{Z}$, and 
 $k$ is an integer such that 
$v^k$ has order $2n+1$. 
Put $X=\{\, \overline{kl} \mid l=1, 2, \ldots, 2n \,\}$, 
where $\overline{kl}=kl \bmod{2n+1}$. 
If $\overline{kl} = \overline{k l^\prime }$ for $l \neq l^\prime$ 
($0<l, l^\prime <2n+1$), 
then $k |l-l^\prime|  \equiv 0 \bmod{2n+1}$, and it follows that 
the order of $w$ is a divisor of $|l-l^\prime|<2n$. 
Then the order is smaller than $2n+1$, which  
is a contradiction. Hence, if $l \neq l^\prime$ ($0<l, l^\prime<2n+1$), 
then $\overline{k l} \neq \overline{k l^\prime}$. 
Since $v^k$ has order $2n+1$, $v^{kl} \neq 1$ for $0<l<2n+1$; 
thus $\overline{0} \not\in X$. 
Thus $X=\{\, \overline{1}, \overline{2}, \ldots, \overline{2n} \,\}$, and hence $X$ contains $\overline{1}$. 
Hence $\overline{k l_0}=\overline{1}$ for some integer $l_0$, and it follows that 
$w^{l_0}=\xi^{-1} v^{k l_0} \xi=\xi^{-1} v \xi$. 
If $w=1$, then $w^{l_0}=\xi^{-1} v \xi=1$; hence $v=1$. 
Conversely, if $v=1$, then $w=1$. 
Thus $G_{4,\,2n+1}^\prime/N_w=G_{4,\,2n+1}^\prime/N_v
=\mathbb{Z}/4 \mathbb{Z}$. 
\qed
\\

Using the results of Theorems \ref{Thm4-6} and \ref{Thm4-9}, we have 
the following theorem. 
%
\begin{thm} \label{0730-1}
For an integer $l>1$, $\mathcal{S}_4(\sigma_1 \sigma_3, \Delta^l)$ 
   is not equivalent to either a spun $T^2$-link, a turned spun $T^2$-link, 
or the split union of 
 spun $T^2$-links and turned spun $T^2$-links. 
\end{thm}
\Proof
Since the link group of 
the spun $T^2$-link or 
the turned spun $T^2$-link of a classical link $L$ 
is isomorphic to the link group of $L$ (\cite{Livingston, Boyle}), it is classical. 
Thus 
the link group of the split union of spun $T^2$-links and turned spun $T^2$-links 
is also classical. 
However, the link group of 
$\mathcal{S}_4(\sigma_1 \sigma_3, \Delta^l)$ ($l>1$) is not 
classical by Theorems \ref{Thm4-6} and \ref{Thm4-9}. 
Thus we have the conclusion. 
\qed

   \section{Ribbon torus-covering $T^2$-links} \label{ribbon}
   
   In this section we show Theorem \ref{Prop2-10}: 
for certain basis $mn$-braids, $\mathcal{S}_{mn}(a,b)$ is ribbon. 
As a corollary, we can see that the torus-covering $T^2$-link of Theorem \ref{0730-1} 
is ribbon (Corollary \ref{0704-1}). 

Let $M$ be a disjoint union of a finite number of handlebodies. 
 The image of $M$ into $\mathbb{R}^4$ by 
 an immersion $\phi$ is called a {\it 3-ribbon} (\cite{Yanagawa})
 if the singularity set consists of 
{\it ribbon singularities}, 
 i.e. the self-intersection of $\phi(M)$ consists of a finite number 
of mutually disjoint 
2-disks, and for each 2-disk $D$, 
the preimage $\phi^{-1}(D)$ consists of a pair of 2-disks 
$D^\prime$, $D^{\prime \prime}$ such that 
$D^\prime \cap D^{\prime \prime}=\emptyset$, $D^\prime \subset \mathrm{Int}M$ and 
$\partial D^{\prime \prime}=D^{\prime \prime} \cap \partial M$. 
An oriented surface link is {\it ribbon} if it bounds a 3-ribbon (see \cite{Yanagawa}). 

\begin{figure}
 \includegraphics*{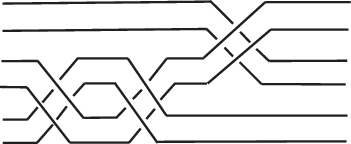}
 \caption{The $2$-parallel $(\sigma_1^2\sigma_2^{-1})^{(2)}$ 
of the 3-braid $\sigma_1^2\sigma_2^{-1}$}
 \label{Nb2}
 \end{figure}

\begin{figure}
 \includegraphics*{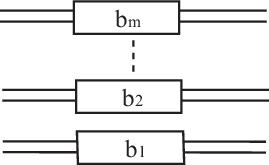}
 \caption{The $mn$-braid $b_1 \circ b_2 \circ \ldots \circ b_m$}
 \label{0728-1}
\end{figure}

For an $m$-braid $b$, we denote by $b^{(n)}$ the $n$-parallel of $b$, i.e. $b^{(n)}$ is 
the $mn$-braid obtained from $b$ by replacing each string of $b$ with its $n$ parallel copies; 
see Fig. \ref{Nb2}. 
 For $n$-braids $b_1,b_2,\ldots, b_m$, 
we denote by $b_1 \circ b_2 \circ \cdots \circ b_m$ the $mn$-braid 
depicted in Fig. \ref{0728-1}. 
 %
\\


   \begin{thm} \label{Prop2-10}
 Let $\alpha$ be a classical $n$-braid whose closure 
$\hat{\alpha}$ is a trivial knot. 
Let $a$ and $b$ be $mn$-braids given by 
\begin{eqnarray*}
&& a=\underbrace{\alpha \circ \alpha 
 \circ \cdots \circ \alpha}_{m}, \\
&& b=b^{\prime (n)} \cdot (\alpha^{l_1} \circ \alpha^{l_2} 
 \circ \cdots \circ \alpha^{l_m}), 
\end{eqnarray*}
where $b^\prime$ is an $m$-braid and $l_j$ is an integer ($j=1,2,\ldots,m$); note that 
$a$ and $b$ are commutative.
 Then $\mathcal{S}_{mn}(a,b)$ is ribbon. 
\end{thm}

\begin{figure}
 \includegraphics*{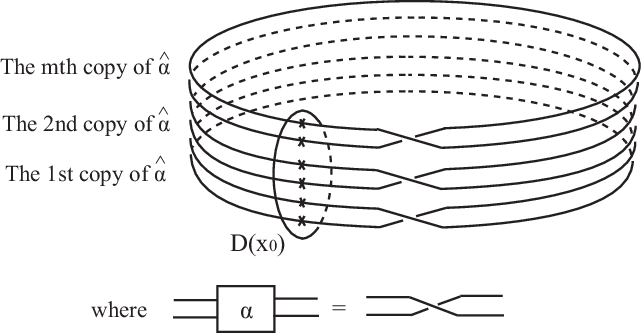}
 \caption{The closure $\hat{a}$ of $a$}
 \label{0831-1}
 \end{figure}

 \Proof 
Let the braid word presentation of $b^\prime$ be 
 $b^\prime=\sigma_{i_1}^{\epsilon_1} \cdot \sigma_{i_2}^{\epsilon_2} \cdots 
 \sigma_{i_{\nu}}^{\epsilon_{\nu}}$, 
  where $i_k \in \{ 1,2,\ldots,m-1 \}$ and $\epsilon_k \in \{+1,-1\}$ for $k=1,2,\ldots,\nu$. 
Let $\pi \,:\, \mathbb{R}^3_+ \times S^1 \rightarrow \mathbb{R}^3_+$ be the projection. 
Let us take a solid torus $N(\mathbf{m}) \subset \mathbb{R}^3_+$ and a disk 
$D(x_0) \subset N(\mathbf{m})$ as 
in the proof of Proposition \ref{eg}. 
Let us take the closure $\hat{a}$ of $a$ in $N(\mathbf{m})$ as shown in Fig. \ref{0831-1}, 
 where we take the $m$ parallel copies of $\hat{\alpha}$ 
in such a position that the identified corresponding ends are in $D(x_0)$. 
 
\begin{figure}
 \includegraphics*{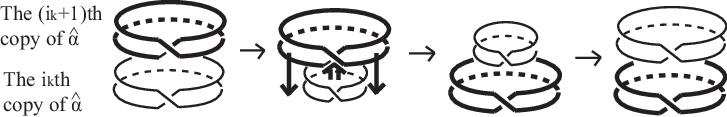}
 \caption{We consider this isotopy (1), 
concerning $(\sigma_{i_k}^{\epsilon_k})^{(n)}$, if $\epsilon_k=+1$.}
 \label{0713-9}
 \end{figure}

\begin{figure}
 \includegraphics*{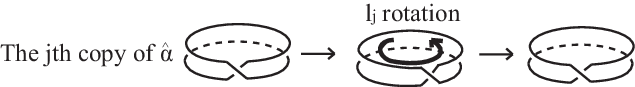}
 \caption{The isotopy (2), concerning 
              $\alpha^{l_1} \circ \alpha^{l_2} \circ \cdots \circ \alpha^{l_m}$}
 \label{0817-1}
 \end{figure}

 We consider a surface link $S$ determined by 
the motion picture $S_t=\pi(S \cap (\mathbb{R}_+^3 \times \{t\}))$ along $S^1$, 
which is the orbit of the isotopy from $\hat{a}$ to $\hat{a}$, 
given by the composition of 
the following (1) and (2). 
\begin{enumerate}[(1)]
\item 
Concerning $b^{\prime (n)}$, 
let us take the isotopy from $\hat{a}$ to $\hat{a}$ as follows. 
For each $(\sigma_{i_k}^{\epsilon_k})^{(n)}$, 
we consider the isotopy shown in Fig. \ref{0713-9} 
if $\epsilon_k=+1$, and its inverse if $\epsilon_k=-1$. 
Further, we consider the composition of them for all $k$. 
\item 
Concerning
$\alpha^{l_1} \circ \alpha^{l_2} \circ \cdots \circ \alpha^{l_m}$, 
let us take the isotopy from $\hat{a}$ to $\hat{a}$ 
which turns the $j$th copy of $\hat{\alpha}$ $l_j$ times as shown in Fig. \ref{0817-1}, 
for each $j=1,2,\ldots,m$. 
\end{enumerate}
Since each isotopy is from $\hat{a}$ to $\hat{a}$, 
$\{S_t\}$, and hence $S$, is well-defined. 
%
%

%
 %
%


\begin{figure}
 \includegraphics*{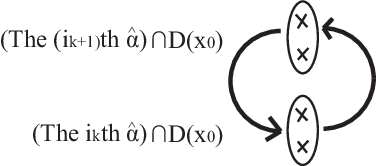}
 \caption{The orbit of $S_t \cap D(x_0)$ by the isotopy (1) 
concerning $(\sigma_{i_k}^{\epsilon_k})^{(n)}$, if $\epsilon_k=+1$} 
 \label{0902-1}
 \end{figure}

Next we show that $S$ is equivalent to $\mathcal{S}_{mn}(a,b)$, as follows. 
It suffices to see that the orbit of $S_t \cap D(x_0)$ forms $b$. 
The orbit of $S_t \cap D(x_0)$ by the isotopy (1) is as in Fig. \ref{0902-1}; thus 
it forms $(\sigma_{i_k}^{\epsilon_k})^{(n)}$.  
Since 
 the the isotopy (2) turns 
the $j$th copy of $\hat{\alpha}$ $l_j$ times ($j=1,2,\ldots, m$), 
by the similar argument to the proof of Proposition \ref{3eg}, 
we can see that the orbit of $S_t \cap D(x_0)$ by this isotopy forms 
$\alpha^{l_1} \circ 
\alpha^{l_2} \circ \cdots \circ \alpha^{l_m}$. 
Thus, the orbit of $S_t \cap D(x_0)$ as a whole forms 
$(\sigma_{i_1}^{\epsilon_1} \sigma_{i_2}^{\epsilon_2} \cdots 
\sigma_{i_\nu}^{\epsilon_\nu})^{(n)} \cdot (\alpha^{l_1} \circ 
\alpha^{l_2} \circ \cdots \circ \alpha^{l_m})$, 
which is $b$. 
 Thus $S$ is equivalent to $\mathcal{S}_{mn}(a,b)$ by Lemma \ref{55} (2). 

\begin{figure}
 \includegraphics*{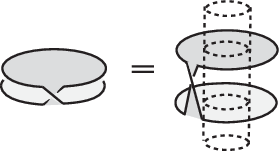}
 \caption{The disk bounded by $\hat{\alpha}$}
 \label{0713-8}
 \end{figure}

Now let us construct an immersed 3-manifold $M$ such that $\partial M=S$, 
which is 
  determined by 
$M_t=\pi(M \cap (\mathbb{R}_+^3 \times \{t\}))$ as follows. 
Since $\hat{\alpha}$ is a trivial knot, 
we can take a disk bounded by $\hat{\alpha}$ 
as shown in Fig. \ref{0713-8}. 
For each $S_t$, let $M_t$ be the union of such disks bounded by $S_t$. 
As the union of $M_t$, we naturally obtain 
 an immersed 3-manifold $M$ such that $\partial M=S$. 

\begin{figure}
 \includegraphics*{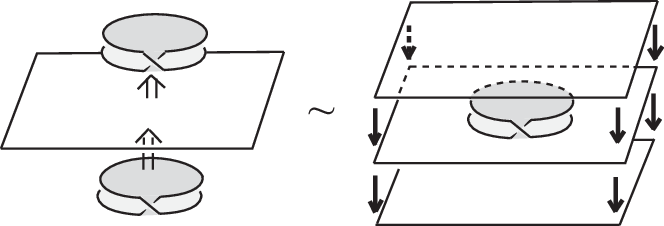}
 \caption{The motion picture of an upward move of the disk 
(the left picture) is equivalent to the motion picture of a downward move of a horizontal plane 
(the right picture). }
 \label{0731-1}
 \end{figure}

In order to show that $\mathcal{S}_{mn}(a,b)$ is ribbon, 
it is sufficient to show that $M$ is a 3-ribbon, i.e. $M$ has only ribbon singularities. 
Since $M$ has no singularity in the motion picture of the isotopy (2), 
it is sufficient to show that 
the singularity in the motion picture of the isotopy (1) 
consists of ribbon singularities. 
Let us consider the singularity of 
$M$ in the motion picture of the isotopy (1) shown in Fig. \ref{0713-9}. 
 This singularity is of the form of the singularity of the motion picture shown in 
Fig. \ref{0731-1}, and hence this singularity set is the disk itself. 
Therefore $M$ has only ribbon sigularities. 
 \qed
 \\

\begin{cor} \label{0704-1}
For any integer $l$, 
 $\mathcal{S}_4(\sigma_1 \sigma_3, \Delta^l)$ 
is ribbon. 
\end{cor}
\Proof
Put $a=\sigma_1 \sigma_3$ and $b=\Delta^l$. 
Let $\alpha$ be a $2$-braid $\sigma_1$. 
Then $\hat{\alpha}$ is a trivial knot, and 
$a=\alpha \circ \alpha$. 
By definition, $\sigma_1^{(2)}=\sigma_2 \sigma_1 \sigma_3 \sigma_2$. 
Since 
$\Delta=(\sigma_2 \sigma_1 \sigma_3 \sigma_2) \cdot \sigma_1 \sigma_3
=\sigma_1^{(2)} \cdot \sigma_1 \sigma_3$, and 
$\sigma_i \cdot \sigma_1^{(2)}=\sigma_1^{(2)} \cdot \sigma_j$ for $\{i,j\}=\{1,3\}$, 
 together with the fact that $\sigma_1$ and $\sigma_3$ commute, 
 it follows that $\Delta^{l}$ can be 
written as 
 $\Delta^{l}=\sigma_1^{l \, (2)} \cdot \sigma_1^l \sigma_3^l$; 
thus $b=b^{\prime (2)} \cdot (\alpha^l \circ \alpha^l)$, 
where $b^\prime=\sigma_1^l$. 
Thus the basis braids $a$ and $b$ 
have the required presentations of Theorem \ref{Prop2-10}, and hence 
$\mathcal{S}_4(\sigma_1 \sigma_3, \Delta^{l})$
 is ribbon by the theorem. 
   \qed
\\

 Thus  
   the torus-covering $T^2$-link of Theorem \ref{0730-1}
   is ribbon. 
Together with Theorem \ref{0730-1}, 
  this indicates the following corollary. 
\begin{cor} \label{0805-1}
 For an integer $l>1$, $\mathcal{S}_4(\sigma_1 \sigma_3, \Delta^l)$ 
has the following properties: 
(1) it can be presented by an $m$-chart on $S^2$ without white vertices, 
however (2) any $m$-chart on $T$ presenting it 
has at least one white vertex. 
\end{cor}

\Proof 
Put $S=\mathcal{S}_4(\sigma_1 \sigma_3, \Delta^l)$. 

(1) Any ribbon surface link is 
presented by an $m$-chart on the standard 2-sphere $S^2$ without white vertices \cite{Kamada92, Kamada3}. 
Since $S$ 
is ribbon by Corollary \ref{0704-1}, it is presented by an $m$-chart on $S^2$ 
without white vertices. 

(2) By Lemma \ref{56} and Proposition \ref{0706-1}, 
if an $m$-chart on $T$ presenting a torus-covering $T^2$-link does not 
have a white vertex, then it presents either 
a spun $T^2$-link, a turned spun $T^2$-link, or  
the split union of spun $T^2$-links and turned spun $T^2$-links. 
Since $S$ 
is not equivalent to such a surface link by Theorem \ref{0730-1}, 
it cannot be presented by an $m$-chart on $T$ without white vertices. 
\qed
\\

Concerning $\mathcal{S}_4( \sigma_1 \sigma_3, \Delta)$, 
we have the following corollary, 
by using the 3-ribbon constructed in the proof of Theorem \ref{Prop2-10}. 

 \begin{cor} \label{0704-6}
The torus-covering $T^2$-knot $\mathcal{S}_4( \sigma_1 \sigma_3, \Delta)$ is unknotted. 
\end{cor}
A {\it 1-handle} attaching to a surface link $S$ is a 3-ball $h$ embedded in $\mathbb{R}^4$ 
such that $S \cap h$ is a pair of 2-disks in $\partial h$. 
The closure (as a set) of 
$(S \cup \partial h)-(S \cap h)$ is a surface link. We call it the surface link obtained 
from $S$ by a {\it 1-handle surgery} along a 1-handle $h$. 
A {\it 2-handle} attaching to $S$ is a 3-ball $h$ embedded in $\mathbb{R}^4$ 
such that $S \cap h$ is an annulus in $\partial h$. The closure (as a set) of 
$(S \cup \partial h)-(S \cap h)$ is a surface link. We call it the surface link obtained 
from $S$ by a {\it 2-handle surgery} along a 2-handle $h$. 
The inverse operation of a 1-handle surgery is a 2-handle surgery, and vice versa. 
\\

\Proof
 By Corollary \ref{0704-1}, 
$\mathcal{S}_4(\sigma_1 \sigma_3, \Delta)$ is 
equivalent to $S=\partial M$ in the proof of Theorem \ref{Prop2-10}. 
We use the notations of the theorem. 
By Corollary \ref{0704-1}, 
$\alpha$ is the $2$-braid $\sigma_1$, and the basis braids are the $4$-braids given by 
$\alpha \circ \alpha$ and 
$(\sigma_1)^{(2)} \cdot (\alpha \circ \alpha)$. 
Since $\partial M_0=S_0$ is the closure of $\alpha \circ \alpha$, 
it consists of two components; thus 
the part of $M$ of the motion picture of the isotopy (1) 
consists of two connected components. 
  Let us denote by $h$ one of 
the components containing the first copy of $\hat{\alpha}$ in $\partial M_0$.  
Since $h$ is an embedded 3-ball 
such that $S \cap h$ is an annulus in $\partial h$ (see Fig. \ref{0713-9}), 
it is a 2-handle on $S$. 
Let $M^\prime=\mathrm{cl}(M-h)$, and put $S^\prime=\partial M^\prime$. 
Then $S^\prime$ is the surface link obtained from $S$ by a 2-handle surgery 
along $h$. 
Since a 2-handle surgery is the inverse operation of a 1-handle surgery, 
$h$ is a 1-handle on $S^\prime$, and 
$S$ is obtained from $S^\prime$ by a 1-handle surgery along $h$. 
Since the singularity set of $M$ is contained in $h$ (see the proof of Theorem \ref{Prop2-10}), 
$M^\prime$ is an embedded 3-ball with no 
singularity; thus $S^\prime$ is an unknotted sphere. 
It is known \cite[Corollary 5]{Boyle88} that if a surface knot is unknotted, then 
the result of a 1-handle surgery for any 1-handle is also unknotted. 
Thus $S$, hence $\mathcal{S}_4(\sigma_1 \sigma_3, \Delta)$, is unknotted. 
\qed

\section{Quandle cocycle invariants} \label{TriplePoint}
 %
%
It is known \cite{Asami-Satoh} that the quandle cocycle invariant
of a twist-spun 2-knot of a classical knot $K$ 
can be presented by using the quandle cocycle invariants of a 1-tangle whose closure is $K$. 
 In this section we present the quandle cocycle invariant 
of $\mathcal{S}_m(b, \Delta^{2n})$ for an $m$-braid $b$ (Theorem \ref{0602-t}), 
using the quandle cocycle invariants of the closure of $b$. 
Here $\Delta$ is a half twist of a bundle of $m$ 
parallel strands. 
In Theorem \ref{Thm2-11}, we calculate some concrete examples. 
 
This section is organized as follows. 
In Section \ref{0901-1}, we review the quandle cocycle invariants and 
the shadow cocycle invariants. 
Further, 
we give a certain 2-cocycle, 
which is determined from a 3-cocycle. 
Using these terms, we give the statement of 
Theorem \ref{0602-t}. 
In Section \ref{0901-2}, we study triple points of $\mathcal{S}_m(b, \Delta^{2n})$,   
and prove Theorem \ref{0602-t}. 
   In Section \ref{s0602-5}, we show Theorem \ref{Thm2-11}, 
using a dihedral quandle and Mochizuki's 3-cocycle. 
  
 \subsection{Quandle cocycle invariant of $\mathcal{S}_m(b, \Delta^{2n})$} \label{0901-1}

  %
  %
  A {\it quandle} (\cite{Joyce}) is a set $Q$ with a binary operation $*$ 
satisfying the following conditions: 
  \begin{enumerate}[(i)]
  \item for any $x \in Q$, $x*x=x$, 
  \item for any $x,\,y \in Q$, there exists a unique $z \in Q$ such that $x=z*y$, and 
  \item for any $x,\,y,\,z \in Q$, $(x*y)*z=(x*z)*(y*z)$. 
  \end{enumerate}
From now on, assume that $Q$ is a finite quandle, i.e. a quandle consisting of finitely many elements. 

\begin{figure}
  \includegraphics*{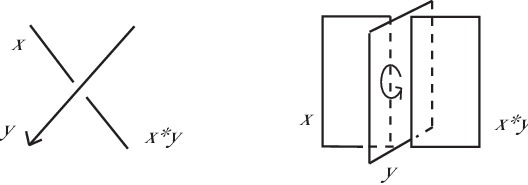}
  \caption{A $Q$-coloring $C$, where $x$, $y$, and $x*y$ are the colors of arcs or broken sheets 
given by $C$}
  \label{1001-1}
  \end{figure}

For an oriented classical link $L$ or 
an oriented surface link $S$, 
let us denote by $D$ the diagram of $L$ or $S$, i.e. 
the image of $L$ or $S$ by a generic projection 
to $\mathbb{R}^2$ or $\mathbb{R}^3$. 
    In order to indicate crossing information of the diagram, 
    we break the under-arc or the under-sheet 
   into two pieces missing the over-arc or the over-sheet. 
   Then the diagram is presented by a disjoint union of arcs, or compact 
  surfaces which are called {\it broken sheets}. 
Let $B(D)$ be the set of such arcs or broken sheets. 
 A {\it $Q$-coloring} for a diagram $D$ of $L$ or $S$ 
is a map $C \,:\, B(D) \rightarrow Q$ 
   as in Fig. \ref{1001-1}. 
The image by $C$ is called the {\it color}. 

Let $G$ be an abelian group. 
A {\it 2-cocycle} with the coefficient group $G$ is a map $f \, :\, Q^2 \rightarrow G$ 
satisfying 
\begin{eqnarray*}
& f(s,u)+f(s*u, t*u)=f(s,t)+f(s*t, u), \ \mathrm{and} & \\
& f(s,s)=0 &
\end{eqnarray*}
for any $s,t,u \in Q$. 
A {\it 3-cocycle} is a map $f \, :\, Q^3 \rightarrow G$ 
satisfying 
\begin{eqnarray*}
& f(s,t,u)+f(s*u, t*u, v)+f(s,u,v)=f(s*t, u,v)+f(s,t,v)+f(s*v, t*v, u*v), \\
& f(s,s,t)=0 \ \mathrm{and} \ f(s, t, t)=0 &
\end{eqnarray*}
for any $s,t,u,v \in Q$. 

\begin{figure}
  \includegraphics*{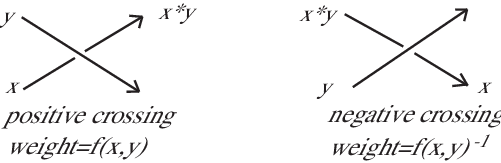}
  \caption{The weight at a crossing, where $x$, $y$, and $x*y$ are the colors of arcs 
by $C$, 
and $f$ is a 2-cocycle}
  \label{0925-1}
  \end{figure}

\begin{figure}
  \includegraphics*{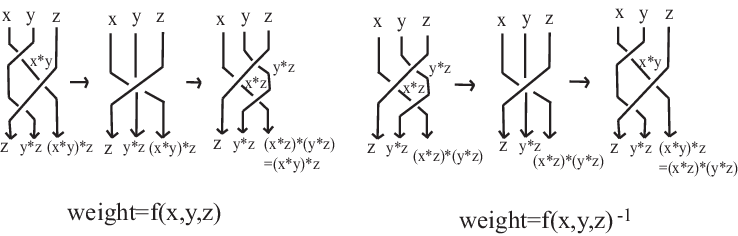}
  \caption{The weight at a triple point, where the triple point is presented by 
a motion picture around it, and $x$, $y$, $z$, etc. are the colors 
by $C$, 
and $f$ is a 3-cocycle}
  \label{0925-2}
  \end{figure}

For a $Q$-coloring $C$ of the diagram $D$ of a classical link $L$ or a surface link $S$, 
we briefly review 
the quandle cocycle invariant as follows (for details see \cite{CJKLS}), 
where $G$ is written multiplicatively. 
For the case of a classical link, 
at each crossing $r$ of the diagram $D$, 
the {\it weight} $W_f(r;C)$ at $r$ for a 2-cocycle $f$ 
is given as in Fig. \ref{0925-1}. 
Put 
\[
\Phi_f(L; C)=\prod_{r \in X_2(D)} W_f(r; C),
\]
where $X_2(D)$ is the set of the crossings of $D$. 
For the case of a surface link, 
at each triple point $t$ of the diagram $D$, 
the {\it weight} $W_f(t;C)$ at $t$ for a 3-cocycle $f$ 
is given as 
 in Fig. \ref{0925-2}   
(\cite[Sections 10 and 11]{CJKLS}, see also \cite[Proposition 4.43 (3)]{Carter-Saito}). 
Put 
\[
\Phi_f(S; C)=\prod_{t \in X_3(D)} W_f(t; C),
\] 
where $X_3(D)$ is the set of the triple points of $D$. 
%
%
%
%
 It is known \cite{CJKLS} that $\Phi_{f}(L;C)$ or $\Phi_f(S;C)$ is an invariant of $L$ or $S$. 
We call it the {\it quandle cocycle invariant} of $L$ or $S$ associated with a $Q$-coloring $C$ 
(see \cite{CJKLS}). 
Since $B(D)$ is a finite set, 
  so is the set of $Q$-colorings for $D$. 
  Let $\mathrm{Col}_Q(D)$ be the set of all the $Q$-colorings. 
Then 
  define $\Phi_{f}(L)$ or $\Phi_f(S)$ by 
  \begin{eqnarray*}
\Phi_f(X) &=& \sum_{C \in \mathrm{Col}_Q(D)} \Phi_f(X;C) \in \mathbb{Z}[G], 
\end{eqnarray*}
where $X=L$ or $S$, and $f$ is a 2-cocycle (resp. 3-cocycle) if $X=L$ (resp. $X=S$). 
  We call $\Phi_f(X)$ the {\it quandle cocycle invariant} of $X$ 
associated with $f$ (\cite{CJKLS}). 
\\

Next we define a shadow coloring of a classical link. 
For a classical link with a given $Q$-coloring, 
its shadow color is determined from the color of the unbounded region (\cite{CKS}), 
which we will call the {\it base color}. 

\begin{figure}
  \includegraphics*{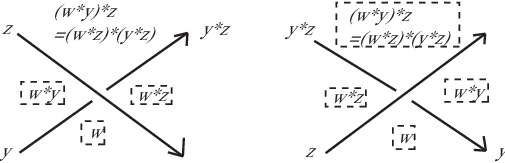}
  \caption{The shadow coloring}
  \label{shadow}
  \end{figure}

  Let $C$ be a $Q$-coloring of the diagram $D$ of a classical link $L$. 
A {\it shadow coloring of $D$ extending $C$ with the base color $x \in Q$} is a map 
$C^*_x \,:\, B^*(D) \rightarrow Q$, where $B^*(D)$ is the union of 
$B(D)$ and the set of regions 
of $\mathbb{R}^2$ separated by the immersed strings of the diagram $D$, satisfying the following conditions: 
\begin{enumerate}[(i)]
\item
$C^*_x$ restricted to $B(D)$ is coincident with $C$. 
\item 
The color of the regions are as in Fig. \ref{shadow}.   
\item 
The color of the unbounded region is $x$. 
\end{enumerate} 
  By \cite{CKS}, $C^*_x$ exists uniquely for given $C$ and $x$. 
For a 3-cocycle $f$, let us 
define the weight at a positive (resp. negative) crossing $r$ by 
$W_f^*(r; C, x)=f(w,y,z)$ (resp. $f(w,y,z)^{-1} \in G$, where 
$y$, $z$, and $w$ are the colors shown in Fig. \ref{shadow}. 
Put 
\[
\Phi_f^*(L;C,x)=\prod_{r \in X_2(D)} W_f^*(r; C, x). 
\]
It is known \cite{CKS} that $\Phi_{f}^*(L;C,x)$ is an invariant of $L$. 
We will call $\Phi_f^*(L; C,x)$ the {\it shadow cocycle invariant} of $L$ 
associated with the $Q$-coloring $C$ and the base color $x$ (see \cite{CKS}). 
 \\

Let $\mathcal{R}_y \,:\, Q \rightarrow Q$ 
be a map defined by $\mathcal{R}_y(x)=x*y$ for $x,y \in Q$. 
Further, let $\mathcal{R}_{\emptyset}=\mathrm{id}_Q$. 
%
We will denote $\mathcal{R}_{y_l} \circ \cdots \circ 
\mathcal{R}_{y_2} \circ \mathcal{R}_{y_1}$ by  
 $\mathcal{R}_{(y_1, y_2, \ldots, y_l)}$ for $(y_1, y_2, \ldots, y_l) \in Q^l$. 
For quandles $Q$ and $Q^\prime$, a map $\phi \,:\, Q \rightarrow Q^\prime$ is called a 
{\it quandle homomorphism} if $\phi(x*y)=\phi(x)*\phi(y)$ for any $x,y \in Q$. 
By the condition (iii) of a quandle, for any $\mathbf{y}=(y_1, y_2, \ldots, y_l)$, $\mathcal{R}_\mathbf{y}$ 
is a quandle homomorphism.

For a $G$-valued 3-cocycle $f$ and $\mathbf{y}=(y_1, \ldots, y_l) \in Q^l$, 
 let $\hat{f}_{\mathbf{y}} \,:\, Q^2 \rightarrow G$ be the map defined by 
\begin{equation} \label{0901-3}
\hat{f}_{\mathbf{y}}(s,t)=\sum_{j=1}^l 
f(\mathcal{R}_{(y_1, \ldots, y_{j-1})}(s), \mathcal{R}_{(y_1, \ldots, y_{j-1})}(t), y_j). 
\end{equation}

\begin{lem}
Assume that 
$\mathcal{R}_{\mathbf{y}}=\mathrm{id}_Q$. Then the map $\hat{f}_{\mathbf{y}}$ is a 2-cocycle. 
\end{lem}
\Proof
We show that $\hat{f}_{\mathbf{y}}$ satisfies the condition of a 2-cocycle, i.e.  
$\hat{f}_{\mathbf{y}} (s,u)+\hat{f}_{\mathbf{y}} (s*u, t*u)=\hat{f}_{\mathbf{y}} (s,t)+\hat{f}_{\mathbf{y}} (s*t, u)$, as follows. 
Put $F=\hat{f}_{\mathbf{y}} (s,u)+\hat{f}_{\mathbf{y}} (s*u, t*u)-\hat{f}_{\mathbf{y}} (s,t)-\hat{f}_{\mathbf{y}} (s*t, u)$. 
Since $f$ is a 3-cocycle, $f$ satisfies  
$f(s,u,v)+f(s*u, t*u, v)-f(s,t,v)-f(s*t, u,v)=
f(s*v, t*v, u*v)-f(s,t,u)$. 
As we mentioned, $\mathcal{R}_{(y_1, \ldots, y_{j-1})}$ is a quandle homomorphism; thus 
we have 
\begin{eqnarray*}
F &=& \sum_{j=1}^l
\{ f(\mathcal{R}_{(y_1, \ldots, y_{j-1})}(s)*y_j, \mathcal{R}_{(y_1, \ldots, y_{j-1})}(t)*y_j, 
\mathcal{R}_{(y_1, \ldots, y_{j-1})}(u)*y_j) \\
 && {} -f(\mathcal{R}_{(y_1, \ldots, y_{j-1})}(s), \mathcal{R}_{(y_1, \ldots, y_{j-1})}(t), 
\mathcal{R}_{(y_1, \ldots, y_{j-1})}(u)) \}. 
\end{eqnarray*}
Since 
$\mathcal{R}_{(y_1, \ldots, y_{j-1})}(s)*y_j
=\mathcal{R}_{(y_1, \ldots, y_j)}(s)$, we have 
\begin{eqnarray*}
 F&=& f(\mathcal{R}_{\mathbf{y}}(s), \mathcal{R}_{\mathbf{y}}(t), 
\mathcal{R}_{\mathbf{y}}(u)) -f(\mathcal{R}_{\emptyset}(s), \mathcal{R}_{\emptyset}(t), 
\mathcal{R}_{\emptyset}(u)), 
\end{eqnarray*}
which is zero from   
  $\mathcal{R}_{\emptyset}=\mathrm{id}$ and the assumption 
$\mathcal{R}_{\mathbf{y}}=\mathrm{id}$. 
Thus $\hat{f}_{\mathbf{y}}$ is a 2-cocycle. 
\qed
\\

For $\mathbf{x}=(x_1, \ldots, x_m)$ and 
$\mathbf{x}^\prime=(x_1^\prime, \ldots, x_m^\prime)$, 
let us denote 
$(x_1, \ldots, x_m, x_1^\prime, \ldots, x_m^\prime)$ 
by $\mathbf{x x^\prime}$. 
\begin{thm} \label{0602-t}
For a given $Q$-coloring of $\hat{b}$, 
let $x_i$ ($i=1,2,\ldots, m$) 
be the color of the $i$th initial arc of the $m$-braid $b$. 
Put $\mathbf{x}=(x_1, \ldots, x_m)$, and 
put 
$\mathbf{y}=\mathbf{x}_1 \mathbf{x}_2 \cdots \mathbf{x}_n$, where 
$\mathbf{x}_1=\mathbf{x}$, and 
$\mathbf{x}_j=\mathcal{R}_{\mathbf{x}_{j-1}}(\mathbf{x}_{j-1})$ ($j>1$). 
Assume that for any $C \in 
\mathrm{Col}_Q(\hat{b})$, $\mathcal{R}_\mathbf{y}=\mathrm{id}_Q$. 
Then the quandle cocycle invariant of $\mathcal{S}_m(b, \Delta^{2n})$ associated with 
a 3-cocycle $f$ is presented by 
\begin{equation*} 
\Phi_f(\mathcal{S}_m(b, \Delta^{2n}))=
\sum_{C \in \mathrm{Col}_{Q}(\hat{b})} 
\Phi_{\hat{f}_\mathbf{y}}(\hat{b}; C) \cdot
 \prod_{i=1}^m \prod_{j=1}^n 
\Psi_{f}^*(\hat{b}; \mathcal{R}_{\mathbf{x}}^{j-1} (C), \mathcal{R}_{\mathbf{x}}^{j-1} (x_i))^{-1}, 
\end{equation*}
where 
$\Phi_{\hat{f}_\mathbf{y}}(\hat{b}; C)$ is the quandle cocycle invariant of $\hat{b}$,  
and $\Psi_{f}^*(\hat{b}; C, x)$ is the shadow cocycle invariant of $\hat{b}$. 
Here 
$\mathbf{x}$ is determined from $C$ and $b$,  
and $\hat{f}_\mathbf{y}$ is the 2-cocycle determined from $f$ and $\mathbf{y}$ by (\ref{0901-3}). 
\end{thm}

\subsection{Proof of Theorem \ref{0602-t}} \label{0901-2}
\subsubsection{Triple points of $\mathcal{S}_m(b, \Delta^{2n})$} \label{s0602-2}
\begin{figure}
  \includegraphics*{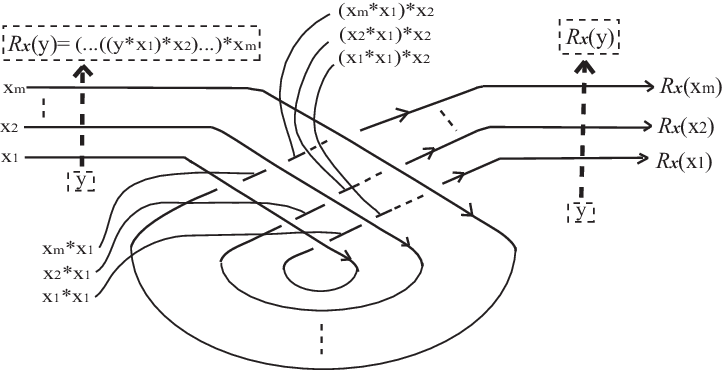}
  \caption{The braid $\Delta^2$}
  \label{0528-2}
  \end{figure}

\begin{figure}
  \includegraphics*{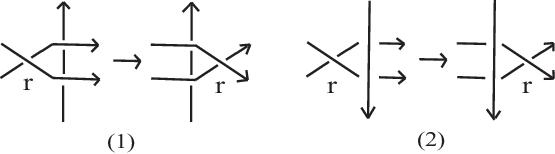}
  \caption{(1) The crossing $r$ slides over an arc, and (2) $r$ slides under an arc. }
  \label{0602-1}
  \end{figure}

Regarding the tubular neighborhood $N(T)$ of 
$T$ as $I \times I \times T$, 
we take for the surface diagram $D$ of $\mathcal{S}_m(b, \Delta^{2n})$ 
the image of the braided surface by the projection 
to $I \times T \subset \mathbb{R}^3 \times \{0\}$. 
Cutting $N(T)$ by $p_T^{-1}(\mathbf{m} \cup \mathbf{l})$, 
we can see that 
$\mathcal{S}_m(b, \Delta^{2n})$ 
is described by a braided surface over a 2-disk 
presenting $b \cdot \Delta^{2n} \cdot b^{-1} \cdot \Delta^{-2n} \rightarrow e$; 
thus $b \cdot \Delta^{2n} \rightarrow \Delta^{2n} \cdot b$, where we use the same notation 
$c$ for the diagram of a classical braid $c$. 
%
 %
Thus the triple points of $D$ appear when we slide $b$ along $\Delta^{2n}$, i.e. 
when we transform $b \cdot \Delta^{2n}$ to $\Delta^{2n} \cdot b$ fixing 
the diagram of $\Delta^{2n}$. 
Each triple point appears when a Reidemeister move of type III occurs. 
Since the braid $\Delta^2$ is isotopic relative the boundary 
to the form as in Fig. \ref{0528-2}, 
 $\mathcal{S}_m(b, \Delta^{2n})$ 
is equivalent to the form such that 
the basis braid $\Delta^{2n}$ is $n$ powers of $\Delta^2$ as in Fig. \ref{0528-2}. 
Since equivalent surface links have the same 
quandle cocycle invariant, we can assume that the diagram of $\Delta^2$ is as in Fig. \ref{0528-2}. 
Sliding $b$ along $\Delta^{2n}$ is equal to 
sliding $b$ along $\Delta^{2}$ $n$ times. 
When we slide $b$ through the $j$th $\Delta^{2}$ 
($j=1,2,\ldots,n$), 
a crossing $r$ of $b$ slides over $m$ arcs, and then under $m$ arcs of $\Delta^{2}$ 
(see Fig. \ref{0602-1}). 
Each time when $r$ slides over or under an arc, a Reidemeister move of type III occurs; thus 
a triple point of $D$ appears. 
Let $t_1^{j,+}(r), \ldots, t_m^{j,+}(r)$, $t_1^{j,-}(r), \ldots, t_m^{j,-}(r)$ 
be the triple points which appear in this order. 


 For a given 3-cocycle $f$ and a $Q$-coloring $C$, we have the following lemma. 
Before a crossing $r$ slides over an arc, 
around $r$ there are three strings. Two strings form $r$, and they separate 
the other string into three arcs. According 
to the orientation, let us call the first arc of the three arcs 
 the {\it initial arc over which $r$ slides}. 
For the colors $x$ and $y$ as in Fig. \ref{0925-1},  
we call the pair $(x, y)$ 
the {\it color of the crossing} $r$ by $C$ (see \cite{CJKLS}). 
\begin{lem} \label{5025}
The weight of $t_i^{j,+}(r)$ is $f(z,x,y)^{-\epsilon}$, and  
the weight of $t_i^{j,-}(r)$ is $f(x,y,z)^{\epsilon}$, 
where $z$ is the color of the initial arc over or under which $r$ slides when it forms 
the triple point, 
and 
$(x,y)$ is the color of $r$ before sliding over or under the arc, and 
$\epsilon=+1$ (resp. $-1$) if $r$ is a positive (resp. negative) crossing.  
\end{lem}

\begin{figure}
  \includegraphics*{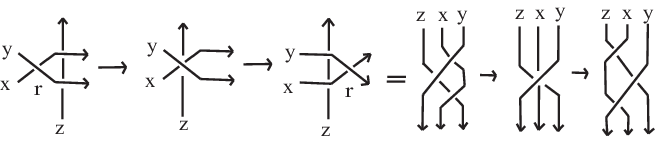}
  \caption{Around $t_+$ if $r$ is a positive crossing}
  \label{0602-4}
  \end{figure}

\begin{figure}
  \includegraphics*{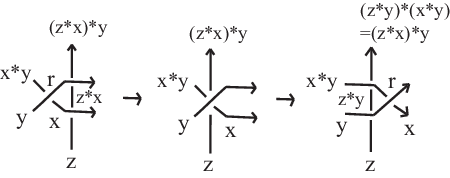}
  \caption{Around $t_+$ if $r$ is a negative crossing}
  \label{1022-1}
  \end{figure}

\begin{figure}
  \includegraphics*{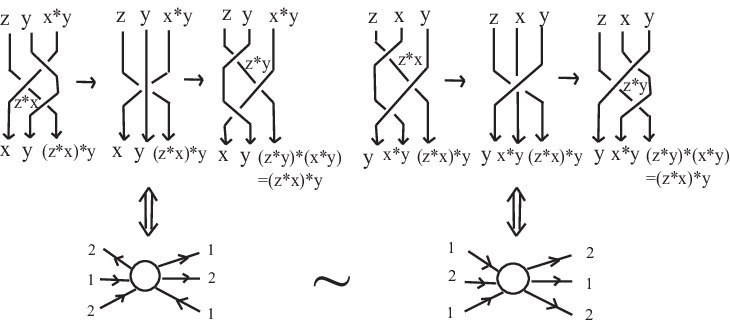}
  \caption{The motion picture of Fig. \ref{1022-1} and its presenting white vertex (the left figure) 
are equivalent to those of the right picture. }
  \label{0602-5}
  \end{figure}

\Proof
Put $t_+=t_i^{j,+}(r)$. 
If $r$ is a positive crossing, then 
the motion picture around $t_+$ is as in Fig. \ref{0602-4}. 
 Thus, for this case, the weight of $t_+$ is $f(z,x,y)^{-1}$; see Fig. \ref{0925-2}. 
If $r$ is a negative crossing, then around $t_+$ is as in 
 Fig. \ref{1022-1}, 
which is equivalent to the right figure of Fig. \ref{0602-5}. 
Thus the weight of $t_+$ is $f(z,x,y)$; see Fig. \ref{0925-2}. 
The weight of $t_i^{j,-}(r)$ is obtained likewise. 
\qed


%
\subsubsection{Proof of Theorem \ref{0602-t}}
We take the surface diagram $D$ of $\mathcal{S}_m(b, \Delta^{2n})$ as in 
Section \ref{s0602-2}. 
First we show that $\mathrm{Col}_Q(D)$ and $\mathrm{Col}_{Q}(\hat{b})$ has one-to-one 
correspondence, as follows. Here, we have the assumption that for 
any $C \in \mathrm{Col}_{Q}(\hat{b})$, $\mathcal{R}_{\mathbf{y}}=\mathrm{id}$. 
For a given $Q$-coloring $C$ of $D$, 
by restricting $C$ to the diagram of the closure of the basis braid $b$, 
we have a $Q$-coloring of $\hat{b}$. 
Conversely, let us consider a given $C \in \mathrm{Col}_{Q}(\hat{b})$. 
Let us give the other basis braid $\Delta^{2n}$ a $Q$-coloring 
such that the colors of the initial arcs are $\mathbf{x}$. 
Since the color of the 
$i$th initial arc of the $j$th $\Delta^{2}$ is the $i$th element of 
$\mathbf{x}_j$ by Lemma \ref{1019-1}, it follows that 
$\mathcal{R}_\mathbf{y}(\mathbf{x})$ are the colors of the terminal arcs 
of $\Delta^{2n}$. Since 
$\mathcal{R}_\mathbf{y}=\mathrm{id}$, 
$C$ can be extended uniquely to the diagram of the closure of $\Delta^{2n}$; thus 
to the closures of 
the basis braids of $\mathcal{S}_m(b, \Delta^{2n})$. 
 Since $\mathcal{S}_m(b, \Delta^{2n})$ is determined from 
the basis braids by Lemma \ref{55} (2), 
$C$ can be extended uniquely to the surface diagram 
$D$. 
 
Now we show the required formula, as follows. 
By definition, 
\[
\Phi_f(\mathcal{S}_m(b, \Delta^{2n}))
=\prod_{r \in X_2(\hat{b})} \prod_{i=1}^m \prod_{j=1}^n \{ W_f(t_i^{j,+}(r)) \cdot W_f(t_i^{j, -}(r)) \}. 
\]

\begin{figure}
  \includegraphics*{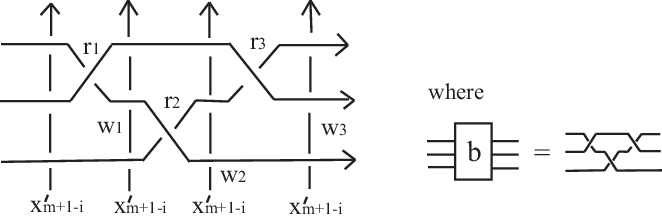}
  \caption{The color $w_k$ of the initial arc over which a crossing $r_k$ slides, 
when it forms $t_i^{j, +}(r_k)$, where $x^\prime_{m+1-i}$ is the color of the $(m+1-i)$th initial arc of the 
$j$th $\Delta^2$. }
  \label{1108-1}
  \end{figure}

First we calculate $\prod_{r \in X_2(\hat{b})} \prod_{i=1}^m \prod_{j=1}^n W_f(t_i^{j,+}(r))$, 
as follows. 
Since the color of the 
$i$th initial arcs of the $j$th $\Delta^{2}$ is $\mathcal{R}_{\mathbf{x}}^{j-1}(x_j)$ 
 ($i=1,2,\ldots,m$, $j=1,2,\ldots,n$) by Lemma \ref{1019-1}, it follows that 
the $Q$-coloring of $b$ before sliding the $j$th $\Delta^2$ is 
$\mathcal{R}_{\mathbf{x}}^{j-1}(C)$. 
The color of a crossing $r \in X_2(\hat{b})$ does not change when 
$r$ slides over an arc. 
When $r$ forms the triple point $t_i^{j,+}(r)$, 
the color of the initial arc over which $r$ slides is the color $w$ depicted in 
Fig. \ref{shadow} (see Fig. \ref{1108-1}), determined from a shadow coloring extending 
$\mathcal{R}_{\mathbf{x}}^{j-1}(C)$ 
with the base color which is the color of the $(m+1-i)$th initial arc of the 
$j$th $\Delta^{2}$, i.e. with the base color $\mathcal{R}_{\mathbf{x}}^{j-1}(x_{m+1-i})$ by 
Lemma \ref{1019-1}. 
Thus $W_f(t_i^{j,+}(r);C)=W_f^*(r; \mathcal{R}_{\mathbf{x}}^{j-1}(C), 
\mathcal{R}_{\mathbf{x}}^{j-1}(x_{m+1-i}))^{-1}$ by Lemma \ref{5025}; 
hence 
$\prod_{r \in X_2(\hat{b})} \prod_{i=1}^m W_f(t_i^{j,+}(r);C) = 
\prod_{i=1}^m \Psi_f^*(\hat{b}; \mathcal{R}_{\mathbf{x}}^{j-1}(C), \mathcal{R}_{\mathbf{x}}^{j-1}(x_i))^{-1}$. 
Hence 
\begin{eqnarray*}
\prod_{r \in X_2(\hat{b})} \prod_{i=1}^m \prod_{j=1}^n W_f(t_i^{j,+}(r);C) 
&=& \prod _{i=1}^m \prod_{j=1}^n 
\Psi^*_f(\hat{b}; \mathcal{R}_{\mathbf{x}}^{j-1}(C), \mathcal{R}_{\mathbf{x}}^{j-1}(x_i))^{-1}. 
\end{eqnarray*}

Next we calculate $\prod_{r \in X_2(\hat{b})} \prod_{i=1}^m \prod_{j=1}^n W_f(t_i^{j,-}(r);C)$, 
as follows. 
For each crossing $r \in X_2(\hat{b})$, the color $(x,y)$ of $r$ by $C$ changes to 
$(x*z, y*z)$, i.e. the color by $\mathcal{R}_z(C)$, 
when 
$r$ slides under an arc as in (2) of Fig. \ref{0602-1}, 
where $z$ is the color of the arc under which $r$ slides. 
Let us denote the $k$th element of $\mathbf{y}=\mathbf{x}_1 \cdots \mathbf{x}_{n}$ by 
$y_k$ ($k=1,2,\ldots, mn$). 
When $r$ slides under the $i$th initial arc 
of the $j$th $\Delta^2$, $r$ has slid under the arcs from the first initial arc 
of the first $\Delta^2$ to the $(i-1)$th initial arc 
of the $j$th $\Delta^2$, whose colors are presented by 
$(y_1, y_2, \ldots, y_{m(j-1)+i-1})$ by Lemma \ref{1019-1}. 
Thus, when $r$ forms the triple point 
$t_i^{j,-}(r)$, the color of $r$ before sliding under an arc is the color by 
$\mathcal{R}_{(y_1, \ldots, y_{m(j-1)+i-1})}(C)$. 
The arc under which $r$ slides is 
the $i$th initial arc of the $j$th $\Delta^2$; thus its color is    
the $i$th element of $\mathbf{x}_j$, i.e. $y_{m(j-1)+i}$, by Lemma \ref{1019-1}. 
Hence  
it follows from Lemma \ref{5025} 
that 
$W_f(t_i^{j,-}(r);C)=
f(\mathcal{R}_{(y_1, \ldots, y_{k-1})}(x), \mathcal{R}_{(y_1, \ldots, y_{k-1})}(y), y_{k})^\epsilon$, 
where $k=m(j-1)+i$ and 
$(x, y)$ is the color of $r$ by $C$, and $\epsilon=+1$ (resp. $-1$) if $r$ is a positive (resp. negative) crossing; 
thus 
$\prod_{i=1}^m \prod_{j=1}^n W_f(t_i^{j,-}(r);C)=
W_{\hat{f}_\mathbf{y}}(r;C)$, and 
we have 
\[
\prod_{r \in X_2(\hat{b})} \prod_{i=1}^m \prod_{j=1}^n W_f(t_i^{j,-};C)=\Phi_{\hat{f}_\mathbf{y}}(\hat{b}; C). 
\]
Hence we have the required formula. 
\qed

  \begin{lem} \label{1019-1} 
In the situation of Theorem \ref{0602-t}, the color of the 
$i$th initial arc of the $j$th $\Delta^{2}$ is the $i$th element of $\mathbf{x}_j$. 
Further, $\mathbf{x}_j=\mathcal{R}_{\mathbf{x}}^{j-1}(\mathbf{x})$.  
\end{lem}

\Proof
The color of the $i$th initial arc of the $j$th $\Delta^{2}$ is the $i$th element of 
$\mathbf{x}$
(resp. $\mathcal{R}_{\mathbf{x}_{j-1}}(\mathbf{x}_{j-1})$) if 
$j=1$ (resp. $j>1$), where $i=1,2,\ldots, m$ and $j=1,2,\ldots, n$, see Fig. \ref{0528-2}; thus 
it is the $i$th element of $\mathbf{x}_j$.  
Since there exits a unique shadow coloring for a given $Q$-coloring and a base color (\cite{CKS}), 
$\mathcal{R}_{\mathcal{R}_\mathbf{x}(\mathbf{x})}(y)=\mathcal{R}_\mathbf{x}(y)$ for any $y \in Q$; 
see Fig. \ref{0528-2}. 
Thus $\mathbf{x}_j=\mathcal{R}_{\mathbf{x}} (\mathbf{x_{j-1}}) 
=\mathcal{R}_{\mathbf{x}}^{j-1}(\mathbf{x})$.  
\qed

\subsection{Concrete calculations} \label{s0602-5}
 The {\it dihedral quandle} of order $p$, denoted by $R_p$, is the 
  set $\{0,1,\ldots,p-1 \}$ with the binary operation $x*y=2y-x \ \pmod{p}$. 
Mochizuki \cite{Mochizuki} showed that for any odd prime $p$, 
the 3-cocycles for $R_p$ with the coefficient group 
$\mathbb{Z}/p \mathbb{Z}$ forms a group isomorphic to $\mathbb{Z}/p \mathbb{Z}$. 
Its generator is reduced (see \cite{Asami-Satoh}) to a map 
given by 
  \[
\theta_p (s,t,u)=v^{(s-t) ( (2u-t)^p+t^p-2u^p )/p} \in 
\langle v \mid v^p=1 \rangle =\mathbb{Z}/p \mathbb{Z}. 
\]
We call $\theta_p$ {\it Mochizuki's 3-cocycle}. 
We identify the group ring $\mathbb{Z}[\mathbb{Z}/p \mathbb{Z}]$ with the Laurent polynomial 
ring $\mathbb{Z}[v, v^{-1}]/(v^p-1)$. 
\begin{thm} \label{Thm2-11}
 We have 
\begin{equation*} 
\Phi_{\theta_p}(\mathcal{S}_4(\sigma_1 \sigma_2^p \sigma_3, \Delta^{2n}))
=p\sum_{i=0}^{p-1} v^{4n i^2} 
   \in \mathbb{Z}[v, v^{-1}]/(v^p=1). 
\end{equation*}
\end{thm}

The {\it triple point number} of a surface link $S$ is the minimum number of triple points 
among all possible diagrams of $S$. 
 By definition, the quandle cocycle invariant of a surface link with the triple point number zero  
has an integer value. 
  Thus we have the following corollary. 
\begin{cor} \label{0526}
  If $n$ is not divisible by $p$, then the triple point number of 
$\mathcal{S}_4(\sigma_1 \sigma_2^p \sigma_3, \Delta^{2n})$ is positive. 
\end{cor}

\noindent
{\it Proof of Theorem \ref{Thm2-11}.} 
 Let us give a $R_p$-coloring for the diagram of 
the basis $4$-braid $b=\sigma_1 \sigma_2^p \sigma_3$. 
By the definition of a $R_p$-coloring, we have $x_1=x_2$ and $x_3=x_4$. 
We will denote the colors by $x$ and $y$ respectively; 
we have $\mathbf{x}=(x,x,y,y)$. 
By a direct calculation, we can see that 
\[
(z*w)*w=z
\]
holds for any $z,w \in R_p$; thus $\mathcal{R}_{(w,w)}=\mathrm{id}$. 
Thus it follows that 
$\mathcal{R}_\mathbf{x}=\mathcal{R}_{(y,y)} \circ \mathcal{R}_{(x,x)}=\mathrm{id}$. 
Hence, by Theorem \ref{0602-t}, 
\begin{equation} \label{0921-1}
\Phi_{\theta_p}(\mathcal{S}_4(\sigma_1 \sigma_2^p \sigma_3, \Delta^{2n}))=
\sum_{C \in \mathrm{Col}_{Q}(\hat{b})}  
 ( \Phi_{\hat{\theta}_{p \mathbf{x}}} (\hat{b}; C) )^n \cdot 
( \Psi_{\hat{\theta}_p}^*(\hat{b}; C,x) \cdot \Psi_{\hat{\theta}_p}^*(\hat{b}; C, y) )^{-2n}, 
\end{equation}
where 
\[
\hat{\theta}_{p \mathbf{x}} (s,t)=\prod_{i=1}^4 
\theta_p(\mathcal{R}_{(x_1, \ldots, x_{i-1})}(s), \mathcal{R}_{(x_1, \ldots, x_{i-1})}(t), x_i).
\]

We calculate $\Phi_{\hat{\theta}_{p \mathbf{x}}}(\hat{b};C)$, as follows. 
Since $x_1=x_2=x$ and $x_3=x_4=y$, we can see that 
$\mathcal{R}_{(x_1, x_2)}=\mathrm{id}$ 
and $\mathcal{R}_{(x_1, x_2, x_3)}=\mathcal{R}_y$. 
Thus $\hat{\theta}_{p \mathbf{x}}(s,t)=\theta_p(s,t, x) \cdot \theta_p(s*x, t*x, x) \cdot 
\theta_p(s,t,y) \cdot \theta_p(s*y, t*y, y) \in \mathbb{Z}/p \mathbb{Z}$. 
By a direct calculation (see \cite{Asami-Satoh}), 
we can see that $\theta_p$ satisfies 
\begin{equation*} \label{0615-1}
\theta_p(s*u, t*u, u)=\theta_p(s,t,u)^{-1}, 
\end{equation*}
for any $s,t,u \in R_p$. 
Thus  
$\hat{\theta}_{p \mathbf{x}}=1$, and hence
$\Phi_{\hat{\theta}_{p \mathbf{x}}} (\hat{b}; C)=1$ for any $C$. 

\begin{figure}
    \includegraphics*{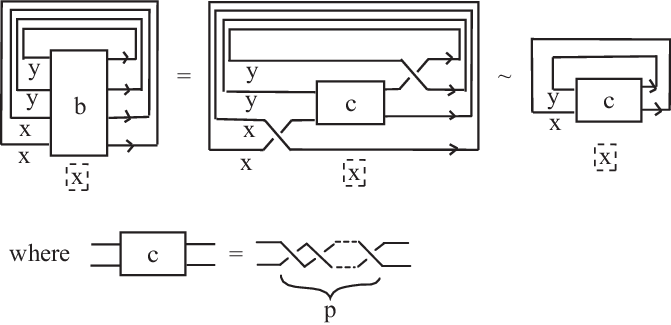}
  \caption{The shadow coloring for $\hat{b}$ with the base color $x$}
  \label{0802-1}
  \end{figure}

\begin{figure}
    \includegraphics*{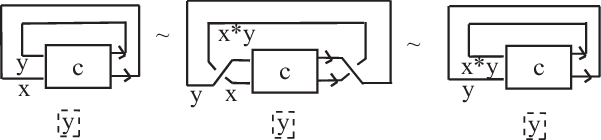}
  \caption{The shadow coloring for $\hat{b}$ with the base color $y$}
  \label{0802-2}
  \end{figure}

We calculate $\Psi_{\hat{\theta}_p}^*(\hat{b}; C,x) \cdot \Psi_{\hat{\theta}_p}^*(\hat{b}; C, y)$, as follows. 
In \cite{Asami-Satoh}, they calculated that $\Psi_{\hat{\theta}_p}^*(\hat{b}; C, x)=v^{-(x-y)^2}$, 
using the diagram of the right figure 
of Fig. \ref{0802-1}. 
Since the diagram of $\hat{b}$ with the $R_p$-coloring $C$ 
with the base color $y$ is transformed as in Fig. \ref{0802-2} 
by Reidemeister moves, $\Psi_{\hat{\theta}_p}^*(\hat{b}; C, y)=v^{-(y-x*y)^2}$, which equals  
$v^{-(y-(2y-x))^2}=v^{-(x-y)^2}$. 
Thus $\Psi_{\hat{\theta}_p}^*(\hat{b}; C,x) \cdot \Psi_{\hat{\theta}_p}^*(\hat{b}; C, y)=v^{-2(x-y)^2}$. 

Hence, by (\ref{0921-1}), 
\begin{eqnarray*} \label{0615-2}
\Phi_{\theta_p}(\mathcal{S}_4(\sigma_1 \sigma_2^p \sigma_3, \Delta^{2n}))
&=& \sum_{x,y \in R_p} v^{4n(x-y)^2} \\
&=& p\sum_{i=0}^{p-1} v^{4n i^2}
   \in \mathbb{Z}[v, v^{-1}]/(v^p=1). 
\end{eqnarray*}
\qed
\\

 The quandle cocycle invariant 
$\Phi_{\theta_p}(\mathcal{S}_4(\sigma_1 \sigma_2^p \sigma_3, \Delta^{2n}))$ has the same value with 
that of the orientation-reversed mirror image of $4n$-twist spun $(2,p)$-torus knot $\tau^{4n}T_p$
(see \cite{Asami-Satoh}). 
An oriented surface link $S$ is {\it invertible} 
if $S$ is equivalent to its orientation-reversed image $-S$, and  
{\it (-)-amphicheiral} if $S$ is equivalent to its orientation-reversed mirror image $-S^*$. 
The $4n$-twist spun 2-knot $\tau^{4n}T_p$
 is equivalent to its mirror image $\tau^{4n}T_p^*$ (see \cite{Litherland}), and in \cite{Asami-Satoh} 
they showed the following fact: 
for an odd prime $p$ with $p \equiv 3 \pmod 4$, if $n$ is not divisible by $p$, 
then $\Phi_{\theta_p}(\tau^{4n}T_p) \neq \Phi_{\theta_p}(-\tau^{4n}T_p^*)$ (see also \cite{CJKS}). 
This means that 
under the above conditions $\tau^{4n}T_p$
 is not invertible. 
Though $\mathcal{S}_4(\sigma_1 \sigma_2^p \sigma_3, \Delta^{2n})$ is invertible for any $p$ and $n$ 
(see Proposition \ref{0803-1}), comparing the quandle cocycle invariants, 
we have the following corollary. 

   \begin{cor} \label{Cor2-12}
   For an odd prime $p$ with $p \equiv 3 \pmod 4$, if $n$ is not divisible by $p$, 
then 
$\mathcal{S}_4(\sigma_1 \sigma_2^p \sigma_3, \Delta^{2n})$ is not (-)-amphicheiral.
   \end{cor}

\acknowledgments
The author would like to thank Professors Takashi Tsuboi and Elmar Vogt for suggesting 
this topic,
and Professors Akio Kawauchi, Tomotada Ohtsuki and the referee for their valuable advice.


\begin{thebibliography}{30}

\bibitem{Artin}
E. Artin, {\it The theory of braids}, Ann. of Math. (2) {\bf 48} (1947), 101--126. 

\bibitem{Asami-Satoh}
S. Asami and S. Satoh, {\it An infinite family of non-invertible surfaces in 4-space}, 
Bull. London Math. Soc. {\bf 37} (2005), 285--296. 

\bibitem{Berstein-Edmonds}
I. Berstein and A. L. Edmonds, {\it On the construction of 
branched coverings of low-dimensional manifolds}, Trans. Amer. Math. Soc. 
{\bf 247} (1979), 87--124. 

\bibitem{B-E84}
I. Berstein and A. L. Edmonds, {\it On the classification of generic 
branched coverings of surfaces}, Illinois J. Math. 
{\bf 28} (1984) 64--82. 
 
 
\bibitem{Birman}
J. Birman, 
{\it Braids, links, and mapping class groups}, 
Annals of Math. Studies, no. 82, Princeton University Press, 1974.

 \bibitem{Bogopoloski}
 O. Bogopolski, {\it Introduction to Group Theory}, European Mathematical Society, 
 2008. 
 
\bibitem{Boyle88}
J. Boyle, {\it Classifying 1-handles attached to knotted surfaces}, Trans. Amer. Math. Soc. {\bf 306} 
(2) (1988), 475--487. 

 \bibitem{Boyle}
 J. Boyle, {\it The turned torus knot in $S^4$}, J. Knot Theory Ramifications 
\textbf{2} (1993), 239--249.

\bibitem{Burde-Zieschang}
G. Burde and H. Zieschang, {\it Eine Kennzeichnung der Torusknoten}, Math. Ann. 
{\bf 169} (1966), 169--176.

\bibitem{Burde-Murasugi}
G. Burde and K. Murasugi,  {\it Links and Seifert fiber spaces}, Duke Math. J. {\bf 37} 
(1970), 89--93.

\bibitem{CJKLS}
J. S. Carter, D. Jelsovsky, S. Kamada, L. Langford and M. Saito, 
{\it Quandle cohomology and state-sum invariants of knotted curves and surfaces}, 
Trans. Amer. Math. Soc.{\bf 355} (2003), 3947--3989.

\bibitem{CKS}
J. S. Carter, S. Kamada, and M. Saito, 
{\it Geometric interpretations of quandle homology and cocycle knot invariant}, 
J. Knot Theory Ramificartions {\bf 10} (2001), 345--358. 


\bibitem{CJKS}
J. S. Carter, D. Jelsovsky, S. Kamada, and M. Saito, 
{\it Computations of quandle cocycle invariants of knotted curves and surfaces}, 
Adv. in Math. {\bf 157} (2001), 36--94.

\bibitem{Carter-Saito}
 J. S. Carter and M. Saito, {\it Knotted surfaces and their diagrams}, Mathematical Surveys 
 and Monographs 55, Amer. Math. Soc., 1998. 
 
\bibitem{Fadell-Neuwirth}
E. Fadell and L. Neuwirth, {\it Configuration spaces}, 
Math. Scand. {\bf 10} (1962), 111--118. 


 %
\bibitem{Gluck}
H. Gluck, {\it The embedding of two-spheres in the four-sphere}, 
Trans. Amer. Math. Soc. {\bf 104} (1962), 308--333. 

 
\bibitem{Hillman}
J. Hillman, {\it 2-Knots and their Groups}, Australian Mathematical Society Lecture Series. 5, 
Cambridge University Press, 1989. 

\bibitem{Hirsch}
M. W. Hirsch, {it Differential Topology}, Graduate Texts in Mathematics 33, 
Springer-Verlag, 1976. 
  
\bibitem{Iwase}
 Z. Iwase, {\it Dehn-surgery along a torus $T^2$-knot}, 
 Pacific J. Math. \textbf{133} (1988), 289--299.
 
\bibitem{Joyce}
D. Joyce, {\it A classical invariant of knots, the knot quandle}, 
J. Pure Appl. Algebra {\bf 23} (1982), 37--65.

\bibitem{Kamada92}
S. Kamada, {\it Surfaces in $R^4$ of braid index three are ribbon}, 
J. Knot Theory Ramifications {\bf 1} (1992), 137--160. 

\bibitem{Kamada92-2}
S. Kamada, {\it 2-dimensional braids and chart descriptions}, 
\lq\lq Topics in Knot Theory (Erzurum, 1992)", 277--287, NATO Adv. Sci. Inst. Ser. 
C Math. Phys. Sci., 399, Kluwer Acad. Publ., (Dordrecht, 1993). 

\bibitem{Kamada1}
S. Kamada, {\it A characterization of groups of closed orientable surfaces in 4-space}, 
Topology {\bf 33} (1994), 113-122.

\bibitem{Kamada2}
S. Kamada, {\it An observation of surface braids via chart description}, 
J. Knot Theory Ramifications {\bf 4} (1996), 517--529.

 \bibitem {Kamada3}
 S. Kamada, {\it Braid and Knot Theory in Dimension Four},
Math. Surveys and Monographs 95, Amer. Math. Soc., 2002. 

\bibitem{Kawauchi}
A. Kawauchi, {\it On pseudo-ribbon surface-links}, J. Knot Theory Ramifications {\bf 11} 
no.7 (2002), 1043--1062.

\bibitem{Litherland}
 R. A. Litherland, {\it Symmetries of twist-spun knots}, in \lq\lq Knot Theory and Manifolds,
Vancouver, B.C., 1983,'' Lecture Notes in Mathematics, {\bf 1144}, pp. 97--107, Springer,
Berlin/New York, 1985. 

\bibitem{Livingston}
C. Livingston, {\it Stably irreducible surfaces in $S^4$}, Pacific J. Math. {\bf 116} 
(1985), 77--84. 

\bibitem{Mochizuki} 
T. Mochizuki, {\it Some calculationsof cohomology groups of finite Alexander quandles}, 
J. Pure Appl. Algebra {\bf 179} (2003), 287--330. 

\bibitem{Montesinos}
J. M. Montesinos, {\it On twins in the four-sphere I}, Quart. J. Math. Oxford (2) {\bf 34} 
(1983), 171--199. 

\bibitem{Neumann}
H. Neumann, {\it Generalized free products with amalgamated subgroups}, 
Amer. J. Math. {\bf 71} (1949), 491--540.

\bibitem{Rudolph}
L. Rudolph, {\it Braided surfaces and Seifert ribbons for closed braids}, 
Comment. Math. Helv. {\bf 58} no.1 (1983), 1--37. 

\bibitem{Satoh-Shima}
S. Satoh and A. Shima, {\it The 2-twist-spun trefoil has the triple 
point number four}, Trans. Amer. Math. Soc. {\bf 356} (2004), 1007-1024. 

\bibitem{Shima}
 A. Shima, {\it Knotted Klein bottles with only double points}, 
 Osaka J. Math. {\bf 40} (2003), 779--799. 
 
 \bibitem{Teragaito}
  M. Teragaito, {\it Symmetry-spun tori in the four sphere}, 
Knots 90, 163--171. 

\bibitem{Yanagawa}
T. Yanagawa, {\it On ribbon 2-knots}, Osaka J. Math. {\bf 6} (1969), 447--464. 

\end{thebibliography}
  \end{document}